\RequirePackage{fix-cm} 
\RequirePackage{fixltx2e}
\RequirePackage[reqno]{amsmath}
\RequirePackage[l2tabu, orthodox]{nag}
\documentclass[final, pdftex, a4paper, twoside, reqno, 12pt]{amsart}
\usepackage{etex}

\usepackage[T1]{fontenc}
\usepackage[latin1]{inputenc}
\usepackage{lmodern}

\usepackage[activate, final, kerning, spacing, factor = 1100, stretch = 10, shrink = 10]{microtype}
\SetTracking{encoding = {*}, shape = sc}{60}

\usepackage[nonewpage, xindy]{imakeidx}
\makeindex[title=Alphabetical index, intoc]

\usepackage{esint}
\usepackage{amsbsy}
\usepackage{amsthm}
\usepackage{dsfont}
\usepackage{empheq}
\usepackage{manfnt}
\usepackage{pifont}
\usepackage{xspace}
\usepackage{amssymb}
\usepackage{upgreek}
\usepackage{stmaryrd}
\usepackage{etoolbox}
\usepackage{extarrows}
\usepackage{textgreek}

\usepackage[nice]{nicefrac}
\usepackage{siunitx}
\usepackage{fourier-orns}

\usepackage{relsize}
\usepackage{textcomp}

\usepackage{mathtools}
\mathtoolsset{centercolon, mathic}

\usepackage[a4paper, margin=1.0in, bottom=1.4in]{geometry}
\usepackage[dvipsnames, table]{xcolor}
\usepackage[most]{tcolorbox}

\definecolor{bckg}{RGB}{20.8, 20.8, 20.8}
\definecolor{oneblue}{rgb}{0.0, 0.0, 0.85}
\definecolor{Lightblue}{RGB}{214, 214, 214}
\definecolor{bluepigment}{rgb}{0.2, 0.2, 0.6}
\definecolor{charcoal}{rgb}{0.21, 0.27, 0.31}
\definecolor{denimblue}{rgb}{0.08, 0.38, 0.74}
\definecolor{darkelectricblue}{rgb}{0.33, 0.41, 0.47}
\definecolor{katyblue}{rgb}{0.129412, 0.137255, 0.63}

\usepackage{booktabs}
\usepackage{array}
\usepackage{ellipsis}
\usepackage{listings}
\usepackage{multicol}

\usepackage{paralist}
\usepackage{csquotes}

\usepackage{graphicx}
\graphicspath{{figs/}}
\DeclareGraphicsExtensions{.pdf, .eps}

\usepackage{tikz}
\usepackage{tikz-cd}
\usepackage{subfigure}
\usepackage{morefloats}
\usepackage{indentfirst}

\usepackage[labelsep=period,
            font=normalsize,
            labelformat=simple,
            labelfont={bf,sf,color=bluepigment},
            justification=raggedright]{caption}

\usepackage[bottom, perpage, multiple, symbol]{footmisc}
\setfnsymbol{wiley}

\newcommand*{\Title}{\textcolor{bluepigment}{Poisson and Symplectic structures}}
\newcommand*{\Longtitle}{Poisson and Symplectic structures, Hamiltonian action, momentum and reduction}
\newcommand*{\Authors}{\textcolor{bluepigment}{V.~Roubtsov and D.~Dutykh}}
\newcommand*{\plogo}{\textcolor{gray}{{\texttt{arXiv.org} / \textsc{hal}}}}
\newcommand*{\Keywords}{Poisson geometry; symplectic geometry; Hamiltonian mechanics; reduction}

\usepackage{url}  
\usepackage[final=true,
           linktoc=all,
           pdftex=true,
           unicode=true,
           bookmarks=true,
           colorlinks=true,
           breaklinks=true,
           urlcolor=katyblue,
           linkcolor=MidnightBlue,
           citecolor=NavyBlue,
           pdffitwindow=false,
           bookmarksopen=false,
           pdfpagemode=UseNone,
           plainpages=false,
           hyperfigures=true,
           pdfstartview={FitV},
           pdftitle={\Longtitle},
           pdfauthor={Volodya Roubtsov and Denys DUTYKH},
           pdfsubject={Poisson geometry},
           pdfcreator={Dr. D},
           pdfproducer={pdfLaTeX},
           pdfkeywords={\Keywords},
           pdfnewwindow=true,
           backref=true,
           pagebackref=true, pdfa]{hyperref}

\usepackage[sort&compress, capitalise, noabbrev]{cleveref}

\usepackage[all]{hypcap}
\usepackage{doi}

\usepackage[sort&compress, comma, square, numbers]{natbib}
\makeatletter
\renewcommand{\@biblabel}[1]{\textbf{[#1]}}
\makeatother

\usepackage[explicit]{titlesec}

\newcommand\invisiblesection[1]{%
  \addcontentsline{toc}{section}{#1}%
  \sectionmark{#1}}

\titleformat{\section}
  {\color{NavyBlue}\Large\sffamily\bfseries}
  {}
  {0em}
  {\hspace*{-0.5em}\colorbox{bckg!5}{\parbox{\dimexpr\linewidth+0\fboxsep\relax}{\centering\Large\thesection. #1}}}
  [\vspace*{0.33em}]

\titleformat{name=\section, numberless}
  {\color{NavyBlue}\Large\sffamily\bfseries}
  {}
  {0.0em}
  {\hspace*{-0.5em}\colorbox{bckg!5}{\parbox{\dimexpr\linewidth+0\fboxsep\relax}{\centering#1}}}
  [\vspace*{0.33em}]

\titleformat{\subsection}
  {\color{NavyBlue}\large\sffamily\bfseries}
  {}
  {0.0em}
  {\hspace*{-0.5em}\colorbox{bckg!5}{\parbox{\dimexpr\linewidth+0\fboxsep\relax}{\centering\thesubsection. #1}}}
  [\vspace*{0.33em}]

\titleformat{name=\subsection, numberless}
  {\color{NavyBlue}\Large\sffamily\bfseries}
  {}
  {0em}
  {\hspace*{-0.5em}\colorbox{bckg!5}{\parbox{\dimexpr\linewidth+0\fboxsep\relax}{\centering#1}}}
  [\vspace*{0.33em}]

\titleformat{\subsubsection}
  {\color{bluepigment}\sffamily\normalsize\bfseries}
  {}
  {0.0em}
  {\hspace*{-0.5em}\colorbox{bckg!5}{\parbox{\dimexpr\linewidth+0\fboxsep\relax}{\centering\thesubsubsection. #1}}}
  [\vspace*{0.33em}]

\titleformat{name=\subsubsection, numberless}
  {\color{bluepigment}\sffamily\normalsize\bfseries}
  {}
  {0.0em}
  {\hspace*{-0.5em}\colorbox{bckg!5}{\parbox{\dimexpr\linewidth+0\fboxsep\relax}{\centering#1}}}
  [\vspace*{0.33em}]

\titleformat{\paragraph}[runin]
  {\color{bluepigment}\sffamily\small\bfseries}
  {}
  {0em}
  {#1}

\titlespacing{\section}{1.0em}{1.5em plus 2pt minus 2pt}%
{1.0em plus 2pt minus 2pt}[0em]
\titlespacing{\subsection}{1.0em}{1.5em plus 2pt minus 2pt}%
{1.0em plus 2pt minus 2pt}[0em]
\titlespacing{\subsubsection}{1.0em}{1.5em plus 2pt minus 2pt}%
{1.0em plus 2pt minus 2pt}[0em]

\usepackage{titletoc}

\contentsmargin{0.5em}
\newlength{\tocsep} 
\titlecontents{section}[\tocsep]
  {\addvspace{10pt}\bfseries\sffamily}
  {\contentslabel[\thecontentslabel]{\tocsep}}
  {}
  {\ \titlerule*[0.75pc]{.}\ \thecontentspage}
  []

\titlecontents{subsection}[\tocsep]
  {\addvspace{8pt}\sffamily}
  {\contentslabel[\thecontentslabel]{\tocsep}}
  {}
  {\ \titlerule*[0.5pc]{.}\ \thecontentspage}
  []

\titlecontents{subsubsection}[\tocsep]
  {\addvspace{4pt}\footnotesize\sffamily}
  {\contentslabel[\thecontentslabel]{\tocsep}\hspace*{0.5em}}
  {}
  {\ \titlerule*[0.35pc]{.}\ \thecontentspage}
  []

\makeatletter
\def\@setauthors{%
  \begingroup
  \def\thanks{\protect\thanks@warning}%
  \trivlist
  \centering\footnotesize \@topsep30\p@\relax
  \advance\@topsep by -\baselineskip
  \item\relax
  \author@andify\authors
  \def\\{\protect\linebreak}%
  \textsc{\normalsize\textcolor{charcoal}{\authors}}%
  \ifx\@empty\contribs
  \else
    ,\penalty-3 \space \@setcontribs
    \@closetoccontribs
  \fi
  \endtrivlist
  \endgroup
}
\def\@settitle{\begin{center}%
  \baselineskip14\p@\relax
    \bfseries
    \textsc{\Large\textcolor{charcoal}{\@title}}
  \end{center}%
}
\makeatother

\usepackage{enumitem}
\setlist[description]{%
  topsep = 9pt,               
  itemsep = 7pt,               
  labelsep = 10pt,
  font={\bfseries\color{NavyBlue}}, 
}
\setlist[itemize]{%
  itemsep = 5pt,
  font={\color{NavyBlue}}
}

\usepackage{fancyhdr}
\usepackage{lastpage}

\numberwithin{equation}{section}

\theoremstyle{plain}
\newtheorem{theorem}{Theorem}[section]
\newtheorem{lemma}{Lemma}[section]
\newtheorem{corollary}{Corollary}[section]
\newtheorem{proposition}{Proposition}[section]

\theoremstyle{definition}
\newtheorem{definition}{Definition}[section]
\newtheorem{example}{Example}[section]

\theoremstyle{remark}
\newtheorem{remark}{Remark}[section]

\usepackage{setspace}
\vfuzz \hfuzz
\newcommand{\up}[1]{$^{\mathrm{\small\textsf{#1}}}$} 

\newtcbox{\mymath}[1][]{%
    nobeforeafter, math upper, tcbox raise base,
    enhanced, colframe = black!35,
    colback = black!5, boxrule = 1pt, arc = 0mm,
    #1}

\newcommand{\comp}{\circ}

\newcommand{\p}{\vec{p}}
\newcommand{\q}{\vec{q}}
\newcommand{\K}{\upkappa}
\newcommand{\A}{\mathds{A}}
\newcommand{\M}{\mathds{M}}

\newcommand{\R}{\mathds{R}}
\newcommand{\T}{\mathds{T}}

\newcommand{\Id}{\mathds{I}}
\newcommand{\ud}{\mathrm{d}}
\newcommand{\ui}{\mathrm{i}}
\newcommand{\ue}{\mathrm{e}}
\newcommand{\vO}{\mathbf{0}}
\renewcommand{\k}{\mathbf{k}}
\renewcommand{\r}{\mathbf{r}}
\renewcommand{\C}{\mathds{C}}
\renewcommand{\L}{\mathbf{L}}
\newcommand{\Ll}{\mathcal{L}}
\renewcommand{\P}{\mathds{P}}
\newcommand{\Rr}{\mathcal{R}}
\renewcommand{\S}{\mathds{S}}
\newcommand{\Xf}{\mathcal{X}}
\newcommand{\Yf}{\mathcal{Y}}
\newcommand{\X}{\mathfrak{X}}
\newcommand{\g}{\mathfrak{g}}
\newcommand{\SU}{\mathbf{SU}}
\renewcommand{\O}{\mathcal{O}}
\newcommand{\omb}{\bar{\omega}}
\newcommand{\so}{\mathfrak{so}}
\newcommand{\su}{\mathfrak{su}}
\newcommand{\dr}{\dot{\mathbf{r}}}

\renewcommand{\mapsto}{\longmapsto}
\newcommand{\ddr}{\ddot{\mathbf{r}}}

\DeclareFontFamily{U}{MnSymbolC}{}
\DeclareSymbolFont{MnSyC}{U}{MnSymbolC}{m}{n}
\DeclareFontShape{U}{MnSymbolC}{m}{n}{
    <-6>  MnSymbolC5
   <6-7>  MnSymbolC6
   <7-8>  MnSymbolC7
   <8-9>  MnSymbolC8
   <9-10> MnSymbolC9
  <10-12> MnSymbolC10
  <12->   MnSymbolC12}{}
\DeclareMathSymbol{\intprod}{\mathbin}{MnSyC}{'270}

\renewcommand{\emptyset}{\varnothing}

\newcommand{\longhookrightarrow}{\lhook\joinrel\longrightarrow}

\DeclareMathOperator{\tr}{tr}

\DeclareMathOperator{\ad}{\textsf{\bfseries ad}}
\DeclareMathOperator{\pr}{\textsf{\bfseries pr}}
\DeclareMathOperator{\Ad}{\textsf{\bfseries Ad}}
\DeclareMathOperator{\GL}{\textsf{\bfseries GL}}
\DeclareMathOperator{\Aut}{\textsf{\bfseries Aut}}
\DeclareMathOperator{\Cas}{\textsf{\bfseries Cas}}
\DeclareMathOperator{\End}{\textsf{\bfseries End}}
\DeclareMathOperator{\Lie}{\textsf{\bfseries Lie}}

\newcommand{\nm}[1]{\textsc{#1}}

\newcommand{\ie}{\emph{i.e.}\xspace}
\newcommand{\eg}{\emph{e.g.}\xspace}

\newcommand{\Mat}{\mathrm{Mat}}

\newcommand{\grad}{\boldsymbol{\nabla}}

\newcommand{\timesb}{\boldsymbol{\times}}

\newcommand{\defeq}{\mathop{\stackrel{\,\mathrm{def}}{\eqcolon}\,}}
\newcommand{\eqdef}{\mathop{\stackrel{\,\mathrm{def}}{\coloneq}\,}}
\newcommand{\pd}[2]{\frac{\partial\hspace{0.0556em} #1}{\partial\/ #2}}

\newcommand{\od}[2]{\frac{\mathrm{d}\hspace{0.0556em} #1}{\mathrm{d}\/#2}}

\newcommand{\odd}[2]{\dfrac{\mathrm{d}\hspace{0.0556em} #1}{\mathrm{d}\/#2}}

\DeclarePairedDelimiterX\abs[1]\lvert\rvert{
  \ifblank{#1}{\:\cdot\:}{\,#1\,}
}
\DeclarePairedDelimiterX\norm[1]\lVert\rVert{
  \ifblank{#1}{\:\cdot\:}{\,#1\,}
}
\DeclarePairedDelimiterX\set[1]{\lbrace}{\rbrace}{\,#1\,}
\DeclarePairedDelimiterX\Set[2]{\lbrace}{\rbrace}{\,#1\ \delimsize\vert\ #2\,}
\DeclarePairedDelimiterX\Inner[2]{\langle}{\rangle}{\,#1\,,\,#2\,}
\DeclarePairedDelimiterX\lb[2]{[}{]}{\,\ifblank{#1}{-}{#1}\,,\,\ifblank{#2}{-}{#2}\,}
\DeclarePairedDelimiterX\pb[2]{\lbrace}{\rbrace}{\,\ifblank{#1}{-}{#1}\,,\,\ifblank{#2}{-}{#2}\,}

\makeatletter
\newcommand{\etabar}{\text{\eta@bar}}
\newcommand{\eta@bar}{%
  \vphantom{$\m@th \eta$}%
  \ooalign{%
    $\m@th \eta$\cr
    \hidewidth\kern.25em\smash{\raisebox{-0.7ex}{$\m@th\mathchar'55$}}\hidewidth\cr}%
}
\makeatother

\makeatletter
\renewcommand*\env@matrix[1][\arraystretch]{%
  \edef\arraystretch{#1}%
  \hskip -\arraycolsep
  \let\@ifnextchar\new@ifnextchar
  \array{*\c@MaxMatrixCols c}}
\makeatother

\makeatletter
\renewenvironment{abstract}{%
    \small\thispagestyle{empty}
    \null\vfil
    {\textcolor{RoyalBlue}{\sc\abstractname.}}
    \quotation
    }
{\endquotation\vfil\null\clearpage}
\makeatother

\newcommand{\TheEnd}{
\bigskip\bigskip
\begin{center}
  \Large
  \decofourleft\hspace*{0.5em}\floweroneleft\hspace*{0.5em}\decoone\hspace*{0.5em}\floweroneright\hspace*{0.5em}\decofourright
\end{center}
\bigskip\bigskip}

\begin{document}

\title[\Title]{\Longtitle}

\author[V.~Roubtsov]{Vladimir Roubtsov\textcolor{denimblue}{$^*$}}
\address{\textcolor{denimblue}{\bf V.~Roubtsov:} LAREMA UMR 6093, CNRS \& Universit\'e d'Angers, 2 Boulevard Lavoisier, 49045 Angers cedex 01, France}
\email{\href{mailto:Vladimir.Roubtsov@univ-angers.fr}{Vladimir.Roubtsov@univ-angers.fr}}
\urladdr{\url{https://math.univ-angers.fr/~volodya/}}
\thanks{\textcolor{denimblue}{$^*$}\it Corresponding author}

\author[D.~Dutykh]{Denys Dutykh}
\address{\textcolor{denimblue}{\bf D.~Dutykh:} Univ. Grenoble Alpes, Univ. Savoie Mont Blanc, CNRS, LAMA, 73000 Chamb\'ery, France and LAMA, UMR 5127 CNRS, Universit\'e Savoie Mont Blanc, Campus Scientifique, 73376 Le Bourget-du-Lac Cedex, France}
\email{\href{mailto:Denys.Dutykh@univ-smb.fr}{Denys.Dutykh@univ-smb.fr}}
\urladdr{\url{http://www.denys-dutykh.com/}}



\begin{titlepage}
\clearpage
\pagenumbering{arabic}
\thispagestyle{empty} 
\noindent
{\Large Vladimir \nm{Roubtsov}}
\\[0.001\textheight]
{\textit{\textcolor{gray}{LAREMA, Universit\'e d'Angers, France}}} 
\\[0.02\textheight]
{\Large Denys \nm{Dutykh}}
\\[0.001\textheight]
{\textit{\textcolor{gray}{CNRS--LAMA, Universit\'e Savoie Mont Blanc, France}}}
\\[0.16\textheight]

\vspace*{2.99cm}

\colorbox{Lightblue}{
  \parbox[t]{1.0\textwidth}{
    \centering\huge
    \vspace*{0.75cm}
    
    \nm{\textcolor{katyblue}{\Longtitle}}
    
    \vspace*{0.75cm}
  }
}

\vfill 

\raggedleft     
{\large \plogo} 
\end{titlepage}


\clearpage
\thispagestyle{empty} 
\par\vspace*{\fill}   
\begin{flushright} 
{\textcolor{RoyalBlue}{\nm{Last modified:}} \today}
\vspace*{1.0em}
\end{flushright}


\clearpage
\maketitle
\thispagestyle{empty}


\begin{abstract}

This manuscript is essentially a collection of lecture notes which were given by the first author at the Summer School \nm{Wis\l{}a}--2019, \nm{Poland} and written down by the second author. As the title suggests, the material covered here includes the \nm{Poisson} and symplectic structures (\nm{Poisson} manifolds, \nm{Poisson} bi-vectors and \nm{Poisson} brackets), group actions and orbits (infinitesimal action, stabilizers and adjoint representations), moment maps, \nm{Poisson} and \nm{Hamiltonian} actions. Finally, the phase space reduction is also discussed. The very last section introduces the \nm{Poisson}--\nm{Lie} structures along with some related notions. This text represents a brief review of a well-known material citing standard references for more details. The exposition is concise, but pedagogical. The Authors believe that it will be useful as an introductory exposition for students interested in this specific topic.

\bigskip\bigskip
\noindent \textcolor{RoyalBlue}{\textbf{\keywordsname:}} \Keywords \\

\smallskip
\noindent \textcolor{RoyalBlue}{\textbf{MSC:}} 
\smallskip \\
\noindent \textcolor{RoyalBlue}{\textbf{PACS:}} 

\end{abstract}


\newpage
\pagestyle{empty}
\tableofcontents
\clearpage
\pagestyle{fancy}


\bigskip\bigskip
\section{Introduction}

These lectures have been delivered by the first author in the Summer School ``\textit{Differential Geometry, Differential Equations, and Mathematical Physics}'' at  \nm{Wis\l{}a}, \nm{Poland} from the 19\up{th} to the 29\up{th} of August 2019. The second author took the notes of these lectures. 

As the title suggests, the material covered here includes the \nm{Poisson} and symplectic structures (\nm{Poisson} manifolds, \nm{Poisson} bi-vectors and \nm{Poisson} brackets), group actions and orbits (infinitesimal action, stabilizers and adjoint representations), moment maps, \nm{Poisson} and \nm{Hamiltonian} actions. Finally, the phase space reduction is also discussed. The \nm{Poisson} structures are a particular instance of \nm{Jacobi} structures introduced by A.~\nm{Lichnerowicz} back in 1977 \cite{Lichnerowicz1977}. Several capital contributions to this field were made by A.~\nm{Weinstein}, see \eg \cite{Weinstein1983}.

The text below does not pretend to provide any new scientific results. However, we believe that this point of view and exposition will be of some interest to our readers. As other general (and excellent) references on this topic include:
\begin{itemize}
  \item R.~\nm{Abraham} and J.E.~\nm{Marsden}. \textit{Foundations of Mechanics}, 2\up{nd} ed., Addison--Wesley Publishing Company, Redwood City, CA, 1987 \cite{Abraham1987}
  \item V.I.~\nm{Arnold}. \textit{Mathematical methods of classical mechanics}. 2\up{nd} ed., Springer, New York, 1997 \cite{Arnold1997}
  \item A.M.~\nm{Vinogradov} and B.A.~\nm{Kupershmidt}. \textit{The structures of Hamiltonian mechanics}, Russ. Math. Surv., \textbf{32}(4), 177--243, 1977 \cite{Vinogradov1977}
\end{itemize}
We can mention also another set of recent lecture notes on the symplectic\index{symplectic geometry} and contact geometries\index{contact geometry} \cite{Tortorella2019}. We mention also the classical review of this topic \cite{Weinstein1998}.

Our presentation remains at quite elementary level of exposition. We restricted deliberately our-selves to the presentation of basic notions and the state-of-the-art as it was in 1990 -- 2000. Nevertheless, we hope that motivated students will be inspired to find more advanced and modern material which is inevitably based on these elementary notes.

In the following text each Section corresponds to a separate lecture. This text is organized as follows. The \nm{Poisson} and symplectic structures are presented in Section~\ref{sec:ps}. The group actions and orbits are introduced in Section~\ref{sec:ga}. The moment map, \nm{Poisson} and \nm{Hamiltonian} actions are described in Section~\ref{sec:mm}. Finally, the manuscript is ended with the description of the phase space reduction in Section~\ref{sec:rp}. The very last Section~\ref{sec:last} is a brief (and essentially incomplete) introduction to \nm{Poisson--Lie}\index{Poisson--Lie structure} structures and some related notions. An excellent account of the last topic can be found in the survey paper by Y.~\nm{Kossmann-Schwarzbach} (1997) \cite{Kosmann-Schwarzbach1997}.


\section{Poisson and symplectic structures}
\label{sec:ps}

\nm{Hamiltonian} systems are usually introduced in the context of the symplectic geometry\index{symplectic geometry} \cite{Fernandes2006, Souriau1997}. However, the use of \nm{Poisson} geometry\index{Poisson geometry} emphasizes the \nm{Lie} algebra\index{Lie algebra} structure, which underlies the \nm{Hamiltonian} mechanics\index{Hamiltonian mechanics}.

\subsection{Poisson manifolds}
\index{Poisson manifold}

Let $\M$ be a smooth manifold\index{smooth manifold} with a bracket\index{bracket}
\begin{equation*}
  \pb*{}{}:\ C^{\infty}\,(\,\M\,)\times C^{\infty}\,(\,\M\,)\ \mapsto\ C^{\infty}\,(\,\M\,)\,,
\end{equation*}
which verifies the following properties:
\begin{description}
  \item[$\bullet$ Bi-linearity\index{bi-linearity}] $\pb*{}{}$ is real-bilinear.
  \item[$\bullet$ Anti-symmetry\index{anti-symmetry}] $\pb*{F}{G}\ =\ -\,\pb*{G}{F}\,$.
  \item[$\bullet$ \nm{Jacobi} identity\index{Jacobi identity}] $\pb*{\pb*{F}{G}}{H}\ +\ \pb*{\pb*{G}{H}}{F}\ +\ \pb*{\pb*{H}{F}}{G}\ =\ 0\,$.
  \item[$\bullet$ \nm{Leibniz} identity\index{Leibniz identity}] $\pb*{FG}{H}\ =\ F\,\pb*{G}{H}\ +\ \pb*{F}{H}\,G\,$.
\end{description}
Then, the bracket $\pb*{}{}$ is a \nm{Poisson} bracket\index{Poisson bracket} and the pair $\bigl(\,\M;\,\pb*{}{}\,\bigr)$ will be called a \nm{Poisson} manifold\index{Poisson manifold}. A \nm{Poisson} algebra\index{Poisson algebra} is defined as the following pair $\bigl(\,C^{\infty}\,(\,\M;\,\R\,);\,\pb*{}{}\,\bigr)\,$. Thanks to the first three properties of the \nm{Poisson} bracket, it is not difficult to see that $\bigl(\,C^{\infty}\,(\,\M;\,\R\,);\,\pb*{}{}\,\bigr)$ is also a \nm{Lie} algebra\index{Lie algebra}. The last property of the \nm{Poisson} bracket\index{Poisson bracket} (\ie the \nm{Leibniz} identity\index{Leibniz identity}) implies that it is also a \emph{derivative}\index{derivative} in each of its arguments.

Let $\bigl(\,\M;\,\pb*{}{}\,\bigr)$ be a \nm{Poisson} manifold\index{Poisson manifold} and $H\ \in\ C^{\infty}\,(\,\M;\,\R\,)\,$, then there exists a unique vector field\index{vector field} $\Xf_{\,H}$ such that
\begin{equation*}
  \Xf_{\,H}(G)\ =\ \pb*{G}{H}\,, \qquad \forall G\ \in\ C^{\infty}\,(\,\M;\,\R\,)\,.
\end{equation*}
The vector field $\Xf_{\,H}$ is called the \nm{Hamiltonian} vector field\index{Hamiltonian vector field} with respect to the \nm{Poisson} structure\index{Poisson structure} with $H$ being the \nm{Hamiltonian} function\index{Hamiltonian function}. Let $\X\,(\,\M\,)$ denote the space of all vector fields\index{vector field} on $\M\,$. Then, the just constructed mapping $C^{\infty}\,(\,\M;\,\R\,)\ \longrightarrow\ \X\,(\,\M\,)$ is a \nm{Lie} algebra morphism\index{Lie algebra morphism}, \ie $\Xf_{\,\pb*{F}{G}}\ =\ \lb{\Xf_{\,F}}{\Xf_{\,G}}\,$.

\begin{definition}
A \nm{Casimir} function\index{Casimir function} on a \nm{Poisson} manifold $\bigl(\,\M;\,\pb*{}{}\,\bigr)$ is a function $F\ \in\ C^{\infty}\,(\,\M;\,\R\,)$ such that for all $G\ \in\ C^{\infty}\,(\,\M;\,\R\,)$ one has
\begin{equation*}
  \pb*{F}{G}\ =\ 0\,, \qquad \forall G\ \in\ C^{\infty}\,(\,\M;\,\R\,)\,.
\end{equation*}
\end{definition}

\subsection{Poisson bi-vector}

If $\bigl(\,\M;\,\pb*{}{}\,\bigr)$ is a \nm{Poisson} manifold\index{Poisson manifold}, then there exists a contravariant anti-symmetric two-tensor\index{anti-symmetric contravariant two-tensor}\index{two-tensor} $\pi\ \in\ \Lambda^{2}\,(\T\M)$ or equivalently
\begin{equation*}
  \pi\,:\ \T^{\ast}\M\,\times\,\T^{\ast}\M\ \longrightarrow\ \R
\end{equation*}
such that
\begin{equation*}
  \Inner*{\pi}{\ud F\wedge\ud G}\,(\,z\,)\ =\ \pi\,(\,z\,)\,\bigl(\,\ud F\,(\,z\,);\,\ud G\,(\,z\,)\,\bigr)\ =\ \pb*{F}{G}\,(\,z\,)\,.
\end{equation*}
In local coordinates $(\,z_{1},\,z_{2},\,\ldots,\,z_{n}\,)$ we have the following expression for the \nm{Poisson} bracket\index{Poisson bracket}:
\begin{equation*}
  \pb*{F}{G}\ =\ \sum_{i,\,j}\,\pi^{ij}\;\pd{F}{z_{i}}\;\pd{G}{z_{j}}\,,
\end{equation*}
where $\pi^{ij}\ \eqdef\ \pb*{z_{i}}{z_{j}}$ are called the elements of the \emph{structure matrix}\index{structure matrix} or the \nm{Poisson} bi-vector\index{Poisson bi-vector} of the underlying \nm{Poisson} structure\index{Poisson structure}. For the vector field\index{vector field} we have the corresponding expression in coordinates: 
\begin{equation*}
  \Xf_{H}\ =\ \sum_{i,\,j}\,\pi^{ij}\;\pd{H}{z_{i}}\;\pd{}{z_{j}}\,, \qquad \text{or} \qquad \Xf_{H}^{j}\ =\ \sum_{i}\,\pi^{ij}\;\pd{H}{z_{i}}\,.
\end{equation*}

\subsubsection{Hamilton map and Jacobi identity}

Let $\pi\ =\ (\,\pi^{ij}\,)$ be a \nm{Poisson} bi-vector\index{Poisson bi-vector} on $\M\,$, then there exists a $C^{\infty}\,(\,\M;\,\R\,)-$linear map\index{linear map} $\pi^{\sharp}:\ \T^{\ast}\M\ \longmapsto\ \T\M$ given by
\begin{equation*}
  \pi^{\sharp}\,(\,\alpha\,)\intprod\beta\ =\ \pb*{\pi}{\alpha\wedge\beta}\,(\,z\,)\ =\ \pi\,(\,z\,)\,\bigl(\,\alpha(\,z\,),\,\beta(\,z\,)\,\bigr)\,,
\end{equation*}
where $\intprod$ denotes the usual interior product\index{interior product} or the substitution of a vector field into a form.

If $\alpha\ =\ \ud f$ for some $f\ \in C^{\infty}\,(\,\M;\,\R\,)\,$, then $\pi^{\sharp}\,(\,\ud H\,)\ =\ \Xf_{H}\,$. It is not difficult to see how the map $\pi^{\sharp}$ acts on the basis elements of co-vectors:
\begin{equation*}
  \pi^{\sharp}\,(\,\ud z_i\,)\ =\ \pb*{z_i}{z_j}\,\pd{}{z_j}\,.
\end{equation*}
Finally, we also have the following \nm{Jacobi} identity:
\begin{equation}\label{eq:1}
  \pi^{il}\;\pd{\pi^{jk}}{z_l}\ +\ \pi^{jl}\;\pd{\pi^{ki}}{z_l}\ +\ \pi^{kl}\;\pd{\pi^{ij}}{z_l}\ =\ 0\,.
\end{equation}

\subsection{Symplectic structures on manifolds}

\begin{definition}
A \emph{symplectic form}\index{symplectic form} on a real manifold $\M$ is a non-degenerate closed $2-$form $\omega\ \in\ \Omega^{\,2}\,(\,\M\,)\ \eqdef\ \Lambda^{\,2}\,(\,\T^{\,\ast}\M\,)\,$. Such a manifold is called a \emph{symplectic manifold}\index{symplectic manifold} and it is denoted by a couple $(\,\M;\,\omega\,)\,$.
\end{definition}
Let $\bigl(\,\M;\,\pb*{}{}\,\bigr)$ be a \nm{Poisson} manifold\index{Poisson manifold} with non-degenerate \nm{Poisson} structure\index{Poisson structure} bi-vector\index{bi-vector} $(\,\pi^{ij}\,)$ and the \nm{Hamiltonian} isomorphism\index{Hamiltonian isomorphism} $\pi^{\sharp}$ such that $\pi^{\sharp}\,(\,\alpha\,)\ =\ \Xf\,$. Then, there is the inverse map\index{inverse map} $(\,\pi^{\sharp}\,)^{\,-1}:\ \T\M\ \longrightarrow\ \T^{\,\ast}\M$ is defined by the following relation:
\begin{equation*}
  \Yf\intprod(\,\pi^{\sharp}\,)^{\,-1}(\,\Xf\,)\ =\ \alpha\,(\,\Yf\,)\,.
\end{equation*}
Moreover, the inverse operator\index{inverse operator} $(\,\pi^{\sharp}\,)^{\,-1}$ defines a $2-$form\index{$2-$form} $\omega_{\pi}$ as follows:
\begin{equation*}
  \omega_{\pi}\,(\,\Xf,\,\Yf\,)\ =\ \Inner*{(\,\pi^{\sharp}\,)^{\,-1}\,(\,\Xf\,)}{\Yf}\,.
\end{equation*}

\subsubsection{Darboux theorem and Hamiltonian vector fields}
\index{Darboux theorem}\index{Hamiltonian vector field}

The following Lemma describes some important properties of the just defined form $\omega_{\pi}$ along with the underlying manifold\index{manifold} $\M\,$:
\begin{lemma}
The real manifold\index{real manifold} $\M$ always has an even dimension. The form $\omega_{\,\pi}$ is a symplectic\index{symplectic form} $2-$form\index{$2-$form} on $\M\,$. The \nm{Jacobi} identity\index{Jacobi identity} \eqref{eq:1} is equivalent to $\ud\,\omega_{\,\pi}\ =\ 0\,$.
\end{lemma}
\begin{proof}
Left to the reader as an exercise.
\end{proof}

\begin{theorem}[G.~\nm{Darboux}]
Let $(\,\M;\,\omega\,)$ be a symplectic manifold\index{symplectic manifold}. There exists a local coordinate system\index{local coordinate system} $(q_{1},\,q_{2},\,\ldots,\,q_{n},\,p_{1},\,p_{2},\,\ldots,\,p_{n})\ \defeq\ (\vec{q},\,\vec{p})$ such that $\omega\ =\ \sum_{i\,=\,1}^{n}\ud\,q_{i}\wedge\ud p_{\,i}\,$. Such coordinates are called \emph{canonical}\index{canonical coordinates} or \nm{Darboux} coordinates\index{Darboux coordinates}.
\end{theorem}

\begin{lemma}
Let $(\,\M;\,\omega\,)$ be a symplectic manifold\index{symplectic manifold} and $H\ \in\ C^{\infty}\,(\M;\,\R)$ the \nm{Hamiltonian} function. Then, there is a unique vector field\index{vector field} $\Xf_{H}$ (\ie the \nm{Hamiltonian} vector field\index{Hamiltonian vector field} associated to the \nm{Hamiltonian} $H$) on $\M$ such that
\begin{equation*}
  \Xf_{H}\intprod\omega\ =\ \ud H\,.
\end{equation*}
The \nm{Hamiltonian} vector field\index{Hamiltonian vector field} $\Xf_{H}$ can be written in the canonical coordinates\index{canonical coordinate} $(\,\vec{q},\,\vec{p}\,)$ on $\M$ as
\begin{equation*}
  \Xf_{H}\ =\ \sum_{i}\,\biggl(\,\pd{H}{p_{\,i}}\;\pd{}{q_{i}}\ -\ \pd{H}{q_{i}}\;\pd{}{p_{\,i}}\,\biggr)\,.
\end{equation*}
The \nm{Poisson} bracket\index{Poisson bracket} in these coordinates looks like
\begin{equation*}
  \pb*{F}{G}\ =\ \Xf_{F}\,(G)\ =\ \sum_{i,\,j}\,\biggl(\,\pd{F}{p_{\,i}}\;\pd{G}{q_{i}}\ -\ \pd{F}{q_{i}}\;\pd{G}{p_{\,i}}\,\biggr)\,.
\end{equation*}
\end{lemma}
\begin{proof}
Left to the reader as an exercise.
\end{proof}

\subsubsection{Example}

Let $\S^{\,2}\ \eqdef\ \Set*{(x,\,y,\,z)\ \in\ \R^{3}}{x^{2}\ +\ y^{2}\ +\ z^{2}\ =\ 1}$ be the $2-$sphere\index{2-sphere}, which can be naturally injected\index{injection} in $\R^{3}\,$: $\imath\,:\ \S^{\,2}\ \xhookrightarrow\ \R^{3}\,$. The $2-$form\index{$2-$form} $\omb\ \in\ \Lambda^{2}\,(\,\R^{3}\,)$ is given by
\begin{equation*}
  \omb\ =\ x\,\ud y\wedge\ud z\ +\ y\,\ud z\wedge\ud x\ +\ z\,\ud x\wedge\ud y
\end{equation*}
and $\omega\ \in\ \Lambda^{2}\,(\,\S^{\,2}\,)$ is defined as $\omega\ =\ \imath^{\ast}\,(\,\omb\,)\,$.

\begin{lemma}
The form $\omega$ gives a symplectic structure\index{symplectic structure} on $\S^{\,2}\,$, \ie $\ud\,\omega\ =\ 0$ and this $2-$form\index{$2-$form} is non-degenerate on $\S^{\,2}\,$.
\end{lemma}

\begin{proof}
First of all, we observe that the closedness\index{closed form} of the $2-$form\index{$2-$form} $\omega$ is a straightforward conclusion in view of
\begin{equation*}
  \ud\,\omega\ =\ \ud\,\bigl(\,\imath^{\ast}(\,\omb\,)\,\bigr)\ =\ \imath^{\ast}\bigl(\,\ud\,\omb\,\bigr)\ =\ 0\,,
\end{equation*}
since $\Lambda^{3}(\,\S^{\,2}\,)\ =\ 0\,$. To check that it is non-degenerate, we make a choice of the following charts\index{chart}:
\begin{align*}
  \phi_{N}:\ \S^{\,2}\setminus\set*{N}\ &\longrightarrow\ \R^{2}\,, \\
  (x,\,y,\,z)\ &\mapsto\ \Bigl(\,\frac{x}{1\ -\ z},\,\frac{y}{1\ -\ z}\,\Bigr)\,,
\end{align*}
\begin{align*}
  \phi_{S}:\ \S^{\,2}\setminus\set*{S}\ &\longrightarrow\ \R^{2}\,, \\
  (x,\,y,\,z)\ &\mapsto\ \Bigl(\,\frac{x}{1\ +\ z},\,\frac{y}{1\ +\ z}\,\Bigr)\,.
\end{align*}
It is not difficult to see that their inverses are given by the following maps:
\begin{equation*}
  \phi_{N}^{\,-1}\,:\ (u,\,v)\ \mapsto\ \Bigl(\,\frac{2\,u}{1\ +\ u^{\,2}\ +\ v^{\,2}},\,\frac{2\,v}{1\ +\ u^{\,2}\ +\ v^{\,2}},\,\frac{u^{\,2}\ +\ v^{\,2}\ -\ 1}{1\ +\ u^{\,2}\ +\ v^{\,2}}\,\Bigr)\,,
\end{equation*}
\begin{equation*}
  \phi_{S}^{\,-1}\,:\ (u,\,v)\ \mapsto\ \Bigl(\,\frac{2\,u}{1\ +\ u^{\,2}\ +\ v^{\,2}},\,\frac{2\,v}{1\ +\ u^{\,2}\ +\ v^{\,2}},\,-\,\frac{u^{\,2}\ +\ v^{\,2}\ -\ 1}{1\ +\ u^{\,2}\ +\ v^{\,2}}\,\Bigr)\,,
\end{equation*}
In both coordinate charts\index{coordinate chart} $(u,\,v)$ induced by $\phi_{N,\,S}$ we obtain:
\begin{multline*}
  \imath^{\ast}\,(\,\omb\,)\ =\ (x\comp\imath)\,\ud\,(y\comp\imath)\wedge\ud\,(z\comp\imath)\ +\ (y\comp\imath)\,\ud\,(z\comp\imath)\wedge\ud\,(x\comp\imath)\ +\ (z\comp\imath)\,\ud\,(x\comp\imath)\wedge\ud\,(y\comp\imath)\ =\\ 
  -\,\frac{4}{1\ +\ u^{\,2}\ +\ v^{\,2}}\;\ud\,u\wedge\ud\,v\ \neq\ 0\,.
\end{multline*}
\end{proof}

\subsection{Co-tangent bundle. Liouville form}
\label{sec:liouville}
\index{co-tangent bundle}\index{Liouville form}

Let $\M$ be a smooth $n-$dimensional manifold\index{smooth manifold} and $\varpi\,:\ \T^{\ast}\M\ \longrightarrow\ \M$ is the projection map\index{projection map} whose differential map\index{differential map} is $T_{\varpi}\,:\ \T\T^{\ast}\M\ \longrightarrow\ \T\,\M\,$.

\begin{definition}
A differential $1-$form\index{$1-$form} $\rho$ on $\T^{\ast}\M\,$, which is defined by $\sigma_{\rho}\,:\ \T^{\ast}\M\ \longrightarrow\ \T^{\ast}\T^{\ast}\M$ as follows
\begin{equation*}
  \Inner*{\sigma_{\rho}}{\Xf}\ =\ \Inner*{\rho}{T_{\varpi_{\rho}}\,(\,\Xf\,)}\,, \qquad \X\ \in\ T_{\varpi_{\rho}}\,(\,\T^{\ast}\M\,)
\end{equation*}
is called the \nm{Liouville} (or \emph{action}) form\index{Liouville form}.
\end{definition}
If $(\vec{q},\,\vec{p})$ is a local coordinate system\index{local coordinate system} on $\T^{\ast}\M\,$, then the form $\rho$ can be written as
\begin{equation*}
  \rho\ =\ \vec{p}\,\ud\,\vec{q}\ =\ \sum_{i=1}^{n}p_{i}\,\ud\,q^{i}\,.
\end{equation*}
The $2-$form\index{$2-$form} $\omega\ \in\ \Omega^{2}\,(\,\T^{\ast}\M\,)\,$, $\omega\ =\ \ud\rho\ =\ \ud\vec{p}\wedge\ud\,\vec{q}\ =\ \sum_{i=1}^{n} \ud p_{i}\wedge\ud\,q^{i}$ is the \emph{canonical} symplectic form\index{symplectic form}.

\subsection{Non-symplectic Poisson structures}
\index{Poisson structure}

Let $\g$ be a real finite-dimensional\index{finite-dimensional Lie algebra} \nm{Lie} algebra\index{Lie algebra} and $\g^{\ast}$ be its dual (as a vector space\index{vector space}). If we suppose that $\dim\g\ =\ n\,$, then $\g^{\ast}$ is isomorphic\index{isomorphism} as a real smooth manifold\index{smooth manifold} to $\R^{n}\,$.

\begin{theorem}
There exists a (non-symplectic) \nm{Poisson} structure\index{Poisson structure} on $\g^{\ast}\,$.
\end{theorem}

\subsection{Poisson brackets on dual of a Lie algebra}
\index{Poisson bracket}

A bracket\index{bracket} can be defined on $C^{\infty}\,(\,\g^{\ast}\,)$ by the identification:
\begin{equation*}
  \g\ \simeq\ \g^{\ast\ast}\ \equiv\ C^{\infty}_{\mathrm{lin}}(\,\g^{\ast}\,)\ \subset\ C^{\infty}\,(\,\g^{\ast}\,)\,, \qquad \Xf\ \mapsto\ F_{\Xf}
\end{equation*}
\begin{equation*}
  F_{\Xf}\,(\xi)\ =\ \Inner*{\xi}{\Xf}\ =\ \xi\,(\,\Xf\,)\,.
\end{equation*}
One has $\pb*{F_{\Xf}}{F_{\Yf}}\ =\ F_{\,\lb{\Xf}{\Yf}}\,$. Let $\set*{e_{i}}_{i\,=\,1}^{n}$ be a base of $\g$ and $\set*{C_{\,ij}^{\,k}}$ be structure constants\index{structure constant} of the \nm{Lie} algebra\index{Lie algebra} $\g\,$, \ie
\begin{equation*}
  \lb{e_{i}}{e_{j}}\ =\ \sum_{k\,=\,1}^{n} C_{\,ij}^{\,k}\,e_{k}\,.
\end{equation*}
Let $\set*{F_{i}}_{i\,=\,1}^{n}$ be a dual basis\index{dual basis} of $\g^{\ast}$ and $\set*{\Xf_{i}}_{i\,=\,1}^{n}$ be the coordinate system\index{coordinate system} of $\g^{\ast}\,$, \ie
\begin{equation*}
  \xi\ =\ \Xf_{i}(\,\xi\,) F_{i}\,, \qquad \Xf_{i}(\,F_{j}\,)\ =\ \delta_{ij}\,, \qquad \Xf_{i}\ =\ F_{e_{i}}\,.
\end{equation*}
Now, it is not difficult to see that
\begin{equation*}
  \pb*{\Xf_{i}}{\Xf_{j}}\ =\ F_{\,\lb{e_{i}}{e_{j}}}\ =\ \sum_{k}\,C_{\,ij}^{\,k}\,F_{e_{k}}\ =\ \sum_{k}\,C_{\,ij}^{\,k}\,\Xf_{k}\,.
\end{equation*}
Finally, we can write the coordinate expression of the \nm{Poisson} bracket\index{Poisson bracket} on the dual of a \nm{Lie} algebra\index{Lie algebra}:
\begin{equation*}
  \pb*{F}{G}\ =\ \sum_{i,\,j,\,k}\,C_{\,ij}^{\,k}\;\pd{F}{\Xf_{i}}\;\pd{G}{\Xf_{j}}\;\Xf_{k}\,.
\end{equation*}

\subsubsection{Definition via gradient operator}
\index{gradient operator}

We define the gradient operator\index{gradient operator}
\begin{equation*}
  \grad\,:\ C^{\infty}\,(\,\g^{\ast}\,)\ \longrightarrow\ C^{\infty}\,(\,\g^{\ast}\,)
\end{equation*}
as follows:
\begin{equation*}
  \Inner*{\eta}{\grad F\,(\,\xi\,)}\ \eqdef\ \od{}{t}\;F\,(\,\xi\ +\ t\,\eta\,)\,\Bigr\vert_{\,t\,=\,0}\,, \qquad \forall F\ \in\ C^{\infty}\,(\,\g^{\ast}\,)\,.
\end{equation*}
In coordinates\index{coordinate} one simply has:
\begin{equation*}
  \grad F\ =\ \pd{F}{\Xf_{i}}\;e_{i}\,.
\end{equation*}
Using the gradient operator\index{gradient operator}, the bracket\index{bracket} is defined as
\begin{equation*}
  \pb*{F}{G}\,(\,\xi\,)\ \eqdef\ \Inner*{\xi}{\lb{\grad\,F\,(\,\xi\,)}{\grad\,G\,(\,\xi\,)}}
\end{equation*}
and in coordinates\index{coordinate} we obtain:
\begin{equation*}
  \pb*{\Xf_{i}}{\Xf_{j}}\,(F_{k})\ =\ \Inner*{F_{k}}{\lb{\grad\,\Xf_{i}\,(F_k)}{\grad\,\Xf_{j}\,(F_k)}}\ =\ \Inner*{F_{k}}{\lb{e_{i}}{e_{j}}}\ =\ C_{\,ij}^{\,k}\,.
\end{equation*}

\subsubsection{Definition via canonical structure on $\T^{\ast}(G)$}

Let $\g$ be a \nm{Lie} algebra\index{Lie algebra}, then there exists a unique (connected and simply connected up to an isomorphism) \nm{Lie} group\index{Lie group} $G$ such that $\T_{e}G\ \simeq\ \g\,$. Let $L_{g}\ \in\ \Aut\,(G)$ be the translation by $g\,$, \ie $\forall\,h\ \in\ G\,$:
\begin{align*}
  L_{g}\,:\ G\ &\stackrel{\simeq}{\longrightarrow}\ G \\
  h\ &\mapsto\ g\,h\,.
\end{align*}
Then, we define
\begin{equation*}
  \lambda_{g}\,(h)\ \eqdef\ \bigl(\T L_{g}\,(h)\bigr)^{\ast}\,:\ \T_{gh}^{\ast}\,(G)\ \longrightarrow\ \T_{h}^{\ast}\,(G)\,,
\end{equation*}
which gives rise to the isomorphism\index{isomorphism} $\lambda_{g}\,:\ \T^{\ast}(G)\ \longrightarrow\ \T^{\ast}(G)$ with the inverse\index{inverse map} $\lambda_{g}^{\,-1}\,$. Define a map\index{map}
\begin{align*}
  \lambda\,:\ G\,\times\,\g^{\ast}\ &\longrightarrow\ \T^{\ast}G \\
  (g,\,\xi)\ &\mapsto\ \lambda_{g}^{\,-1}\,(g)\,(\xi)\,,
\end{align*}
which is a diffeomorphism\index{diffeomorphism} in the following commutative diagram\index{commutative diagram}:
\begin{equation*}
\begin{tikzcd}
  G\,\times\,\g^{\ast} \arrow[rr, "\lambda"] \arrow[dr, "\pr_{1}"'] & & \T^{\ast}G \arrow[dl, "\pi"] \\
  & G &
\end{tikzcd}
\end{equation*}
The co-tangent bundle\index{co-tangent bundle} $\T^{\ast}G\ \simeq\ G\,\times\g^{\ast}$ is a trivial vector bundle\index{trivial vector bundle} with the fiber\index{fiber} $\g^{\ast}\,$. The \nm{Liouville} form\index{Liouville form} $\rho\ \in\ \Omega^{1}\,(\T^{\ast}G)$ defines by the section $\sigma_{\rho}\,:\ \T^{\ast}G\ \longrightarrow\ \T^{\ast}\T^{\ast}G$ similar to Section~\ref{sec:liouville}:
\begin{equation*}
  \Inner*{\sigma_{\rho}(\,\xi\,)}{\Xf}\ =\ \Inner*{\rho}{\T_{\pi_{\xi}}\,(\,\Xf\,)}\,, \qquad \Xf\ \in\ \T_{\pi_{\xi}}\,(\,\T^{\ast}G\,)\,, \qquad \xi\ \in\ \T^{\ast}_{h}G\,.
\end{equation*}
Let $g\ \in\ G\,$, $\xi\ \in\ \T^{\ast}_{gh}G\,$, $\Xf\ \in\ \T_{\pi_{\lambda_{g}(\,\xi\,)}}(\,\T^{\ast}G\,)$\,, then we have
\begin{multline*}
  \Inner*{\sigma_{\rho}\comp\lambda_{g}(\,\xi\,)}{\Xf}\ =\ \Inner*{\lambda_{g}(\,\xi\,)}{\T_{\pi_{\lambda_{g}(\,\xi\,)}}}\ =\ \Inner*{(\T_{h}L_{g})^{\ast}\xi}{\T_{\lambda_{g(\,\xi\,)}}}\ =\\ 
  \Inner*{\xi}{(\T_{h}L_{g})\T_{\lambda_{g}(\,\xi\,)}\Xf}\ =\ \Inner*{\xi}{\T_{\lambda_{g}(\,\xi\,)}(L_{g}\comp\pi)\Xf}\,.
\end{multline*}
We can make two observations:
\begin{itemize}
  \item $(L_{g}\comp\pi)\,(h,\,\xi)\ =\ gh\,$,
  \item $\pi\comp\lambda_{g^{-1}}\ =\ gh\,$.
\end{itemize}
Henceforth, $L_{g}\comp\pi\ \equiv\ \pi\comp\lambda_{g^{-1}}\,$. The \nm{Liouville} form\index{Liouville form} becomes:
\begin{multline*}
  \Inner*{\sigma_{\rho}\comp\lambda_{g}(\,\xi\,)}{\Xf}\ =\ \Inner*{\xi}{\T_{\lambda_{g}(\,\xi\,)}(\,\pi\comp\lambda_{g^{-1}}\,)\Xf}\ =\ \Inner*{\xi}{\T_{\xi}\pi\comp\T_{\lambda_{g}(\,\xi\,)\lambda_{g^{-1}}(\,\Xf\,)}}\ =\\ 
  \Inner*{\rho(\,\xi\,)}{\T_{\lambda_{g}(\,\xi\,)}\lambda_{g^{-1}}(\,\xi\,)\Xf}\ =\ \Inner*{\lambda_{g^{-1}}^{\ast}\comp\rho(\,\xi\,)}{\Xf}\,.
\end{multline*}
We have $\omega\ =\ \ud\,\rho$ as the canonical symplectic form on $\T^{\ast}G\,$. We observe also that
\begin{itemize}
  \item $\sigma_{\rho}\comp\lambda_{g}\ =\ \lambda_{g^{-1}}^{\ast}(\,\sigma_{\rho}\,)\,$,
  \item $\omega\comp\lambda_{g}\ =\ \lambda_{g^{-1}}^{\ast}(\,\omega\,)\,$.
\end{itemize}
We recall that a \nm{Poisson} bracket\index{Poisson bracket} $\pb*{F}{G}$ for $F\,$, $G\ \in\ C^{\infty}\,(\,\T^{\ast}G\,)$ can be defined via $\pb*{F}{G}\ \eqdef\ \Xf_{F}(\,G\,)\,$, where $\Xf_{F}$ is the unique vector field on $\T^{\ast}G$ such that $\Xf_{F}\intprod\omega\ =\ \ud F\,$.

\begin{lemma}\label{lem:1}
Let $g\ \in\ G$ and $F\,$, $G\ \in\ C^{\infty}\,(\,\T^{\ast}G\,)\,$, then
\begin{equation*}
  \pb*{F\comp\lambda_{g}}{G\comp\lambda_{g}}\ =\ \pb*{F}{G}\comp\lambda_{g}\,.
\end{equation*}
\end{lemma}
\begin{proof}
Left to the reader as an exercise.
\end{proof}

Let $C^{\infty}\,(\,\T^{\ast}G\,)^{G}$ denote a subspace\index{subspace} of stable or invariant functions with respect to the mapping $\lambda_{g}\,$, \ie
\begin{equation*}
  C^{\infty}\,(\,\T^{\ast}G\,)^{G}\ \eqdef\ \Set*{F\ \in\ C^{\infty}\,(\,\T^{\ast}G\,)}{F\comp\lambda_{g}\ =\ F}\,, \qquad g\ \in\ G\,.
\end{equation*}
Lemma~\ref{lem:1} shows that the set $C^{\infty}\,(\,\T^{\ast}G\,)^{G}$ is closed with respect to the \nm{Poisson} bracket\index{Poisson bracket}.

Let the linear mapping $\Phi$ be defined as
\begin{align*}
  \Phi\,:\ C^{\infty}\,(\,\g^{\ast}\,)\ &\longrightarrow\ C^{\infty}\,(\,\T^{\ast}G\,) \\
  F\ &\mapsto\ F\comp\pr_{2}\,,
\end{align*}
where $\pr_{2}\,: \T^{\ast}G\ \equiv\ G\,\times\,\g^{\ast}\ \longrightarrow\ \g^{\ast}$ is the canonical projection\index{projection} on the second argument. Let $\imath\,:\ \g^{\ast}\ \longrightarrow\ \T^{\ast}G$ be the canonical embedding\index{canonical embedding}. Then,
\begin{align*}
  \Psi\,:\ C^{\infty}\,(\,\T^{\ast}G\,)^{G}\ &\longrightarrow\ \g^{\ast} \\
  F\ &\mapsto\ F\comp\imath
\end{align*}
is linear and inverse\index{inverse map} to $\Phi\,$, \ie $\Psi\ \equiv\ \Phi^{\,-1}\,$.

\begin{lemma}
The bracket\index{bracket} $\pb*{F}{G}\ \eqdef\ \Phi^{\,-1}\bigl(\,\pb*{\Phi(F)}{\Phi(G)}\,\bigr)$ is a \nm{Poisson} bracket\index{Poisson bracket} on $C^{\infty}\,(\,\g^{\ast}\,)$ coinciding with two previous definitions.
\end{lemma}
\begin{proof}
Left to the reader as an exercise.
\end{proof}


\section{Group actions and orbits}
\label{sec:ga}

Let $G$ be a \nm{Lie} group\index{Lie group} and $\M$ is a smooth manifold\index{smooth manifold}.
\begin{definition}
A \emph{left action}\index{left action} of $G$ on $\M$ is a smooth map\index{smooth map} $\mu\,:\ G\times\M\ \longrightarrow\ \M$ such that
\begin{itemize}
  \item $\mu\,(\,e,\,m\,)\ =\ m\,$, $\forall\,m\ \in\ \M\,$,
  \item $\mu\,\bigl(\,g,\,\mu\,(\,h,\,m\,)\bigr)\ =\ \mu\,(\,gh,\,m\,)\,$, $\forall\,g,\,h\ \in\ G$ and $\forall\,m\ \in\ \M\,$.
\end{itemize}
\end{definition}
\begin{definition}
A \emph{right action}\index{right action} of $G$ on $\M$ is a smooth map\index{smooth map} $\rho\,:\ \M\times G\ \longrightarrow\ \M$ such that
\begin{itemize}
  \item $\rho\,(\,m,\,e\,)\ =\ m\,$, $\forall\,m\ \in\ \M\,$,
  \item $\rho\,\bigl(\,\rho\,(\,m,\,h\,),\,g\,\bigr)\ =\ \rho\,(\,hg,\,m\,)\,$, $\forall\,g,\,h\ \in\ G$ and $\forall\,m\ \in\ \M\,$.
\end{itemize}
\end{definition}
Left and right actions\index{group action} of $G$ and $\M$ are in one-to-one correspondence\index{one-to-one correspondence} by the following relation:
\begin{equation*}
  \rho\,(\,m,\,g^{-1}\,)\ =\ \mu\,(\,g,\,m\,)\,.
\end{equation*}
From now on we shall denote the left\index{left action} \nm{Lie} group\index{Lie group} action\index{group action} $\mu\,(\,g,\,m\,)$ simply by $g\cdot m\,$. We can define several important action types:
\begin{description}
  \item[Effective or faithful] $\forall\,g\ \in\ G\,$, $g\ \neq\ e$ $\quad \Longrightarrow \quad \exists\,m\ \in\ \M$ such that $g\cdot m\ \neq\ m\,$.
  \item[Free]  If $g$ is a group element and $\exists\,m\ \in\ \M$ such that $g\cdot m = m$ (that is, if $g$ has at least one fixed point),  $\Longrightarrow g=e$. Note that a free action on a non-empty $M$ is faithful\index{faithful action}\index{effective action}.
  \item[Transitive] If $\forall\,m,\,n\ \in\ \M\,$, $\exists\,g\ \in\ G$ such that $g\cdot m\ =\ n\,$. In this case the smooth manifold\index{smooth manifold} $\M$ is called homogeneous\index{homogeneous manifold}\index{transitive action}.
\end{description}
Important examples of group actions\index{group action} include:
\begin{example}
$G$ acts on itself by left multiplication\index{multiplication}:
\begin{align*}
  G\times G\ &\longrightarrow\ G \\
  (g,\,h)\ &\mapsto\ gh\,.
\end{align*}
This action\index{group action} is effective\index{effective action} and transitive\index{transitive action}. Indeed, $g\cdot h\ =\ h$ $\Longrightarrow$ $g\ =\ e$ and if $g\cdot m\ =\ n$ $\Longrightarrow$ $g\ =\ n\cdot m^{-1}\,$.
\end{example}
\begin{example}
$G$ acts on itself by conjugation:
\begin{align*}
  G\times G\ &\longrightarrow\ G \\
  (g,\,h)\ &\mapsto\ g\cdot h\cdot g^{-1}\,.
\end{align*}
Generally, this action\index{group action} is not free\index{free action}, transitive\index{transitive action} or effective\index{effective action}.
\end{example}
\begin{example}
$\GL_{\,n}(\,\R\,)$ acts on $\R^{\,n}\setminus\set*{\vO}$ by matrix multiplication\index{matrix multiplication} on the left:
\begin{align*}
  \GL_{\,n}(\,\R\,)\times \R^{\,n}\setminus\set*{\vO}\ &\longrightarrow\ \R^{\,n}\setminus\set*{\vO} \\
  (A,\,x)\ &\mapsto\ A\,x\,.
\end{align*}
This is an example of an effective\index{effective action} and transitive action\index{transitive action}.
\end{example}

\subsection{Stabilizers and orbits}

Let $G$ be a \nm{Lie} group\index{Lie group} which acts on a smooth manifold\index{smooth manifold} $\M\,$. The \emph{orbit}\index{orbit} of a point $m\ \in\ \M$ is
\begin{equation*}
  G\cdot m\ \eqdef\ \Set*{g\ \in\ G}{g\cdot m}\ \subseteq\ \M\,.
\end{equation*}
A \emph{stabilizer}\index{stabilizer} of a point $m\ \in\ \M$ is
\begin{equation*}
  G_{m}\ \eqdef\ \Set*{g\ \in\ G}{g\cdot m\ =\ m}\ \subseteq\ G\,.
\end{equation*}

\begin{proposition}
The stabilizer $G_{m}$ is a closed subgroup of $G$ and $G_{g\cdot m}\ =\ g\cdot G_{m}\cdot g^{-1}\,$, $\forall\,g\ \in\ G\,$.
\end{proposition}
\begin{proof}
Left to the reader as an exercise.
\end{proof}

We mention here two technical theorems regarding the orbits\index{orbit} and stabilizers\index{stabilizer}:
\begin{theorem}
Let $G$ be a \nm{Lie} group\index{Lie group} which acts on a smooth manifold\index{smooth manifold} $\M$ and $m\ \in\ \M\,$. There is a manifold structure\index{manifold structure} on the orbit $G\cdot m$ such that the map
\begin{align*}
  G\ &\longrightarrow\ G\cdot m \\
  g\ &\mapsto\ g\cdot m
\end{align*}
is a submersion\index{submersion} and the embedding\index{embedding} $\imath\,:\ G\cdot m\ \hookrightarrow\ \M$ is an immersion\index{immersion}.
\end{theorem}

\begin{theorem}
The \nm{Lie} algebra\index{Lie algebra} $\g_{\,m}$ of the stabilizer $G_{m}$ for a point $m\ \in\ \M$ coincides with $\ker\T_{e}\Phi\,$, where the mapping $\Phi$ is defined as
\begin{align*}
  \Phi\,:\ G\ &\longrightarrow\ \M \\
  g\ &\mapsto\ g\cdot m\,.
\end{align*}
\end{theorem}

\subsection{Infinitesimal action}

Let $\mu\,:\ G\times\M\ \longrightarrow\ \M$ be a \nm{Lie} group\index{Lie group} action\index{group action} on $\M$ and $\g\ =\ \Lie\,(G)$ be its \nm{Lie} algebra\index{Lie algebra}.

\begin{definition}
Let $\Xf\ \in\ \g$ and $\phi\,:\ \R\ \longrightarrow\ G$ its exponential flow\index{exponential flow}, \ie $\phi\,(t)\ =\ \exp(\,t\,\Xf)\,$. Then, there exists the unique vector field\index{vector field} $\Xf_{\,\M}\ \in\ \X\,(\,\M\,)$ with the flow $\phi_{m}\,:\ \R\ \longrightarrow\ \M$ defined by $\phi_{m}\,(t)\ =\ \phi\,(t)\cdot m\,$. The vector field\index{vector field} $\Xf_{\,\M}$ is defined by
\begin{equation*}
  \Xf_{\,\M}(\,m\,)\,(\,f\,)\ \eqdef\ \od{}{t}\;\bigl(\,f\comp\phi(\,t\,)\cdot m\,\bigr)\,\Bigl\vert_{\,t\,=\,0}\,.
\end{equation*}
The mapping\index{mapping} $\mu_{\ast}\,:\ \g\ \longrightarrow\ \X\,(\,\M\,)$ is called the \emph{infinitesimal action}\index{infinitesimal action} of $\g$ on $\M\,$.
\end{definition}

\begin{proposition}
The mapping\index{mapping} $\mu_{\ast}\,:\ \g\ \longrightarrow\ \X\,(\,\M\,)$ is a \nm{Lie} algebra\index{Lie algebra} (anti-)homomorphism\index{homomorphism}\index{anti-homomorphism} (and it is in particular a linear mapping\index{linear mapping}):
\begin{equation*}
  \mu_{\ast}\bigl(\,\lb{\Xf}{\Yf}\,\bigr)\ =\ -\lb{\mu_{\ast}(\,\Xf\,)}{\mu_{\ast}(\,\Yf\,)}\,.
\end{equation*}
\end{proposition}
\begin{proof}
Left to the reader as an exercise.
\end{proof}

\begin{remark}
One can see that $\mu_{\ast}(\,\Xf\,)_{m}\,(\,f\,)\ =\ \Xf\,(f\comp\Phi_{m})\,$. In other words, $\mu_{\ast}(\,\Xf\,)_{m}\ =\ \T_{e}\Phi_{m}(\,\Xf\,)\,$.
\end{remark}

\begin{proposition}
Let $m\ \in\ \M\,$, then
\begin{equation*}
  \T_{m}G\cdot m\ =\ \Set*{\Xf\ \in\ \g}{\mu_{\ast}(\,\Xf\,)_{m}}\,.
\end{equation*}
\end{proposition}
\begin{proof}
Left to the reader as an exercise.
\end{proof}

The following difficult result is left without the proof:
\begin{theorem}[R.~\nm{Palais}]\label{thm:palais}
Let $G$ be a simply connected\index{simply connected group} \nm{Lie} group\index{Lie group} with the \nm{Lie} algebra\index{Lie algebra} $\g\ =\ \Lie\,(G)$ and $\M$ be a smooth compact manifold\index{smooth manifold} such that there exists a homomorphism\index{homomorphism} of \nm{Lie} algebras\index{Lie algebra} $\rho\,:\ \g\ \longrightarrow\ \X\,(\,\M\,)\,$. Then, there is a unique action\index{group action} $\mu\,:\ G\times\M\ \longrightarrow\ \M$ such that $\mu_{\ast}\ =\ \rho\,$.
\end{theorem}

\begin{proposition}
Let $\mu\,:\ G\times\M\ \longrightarrow\ \M$ be an action\index{group action} of $G$ on a smooth manifold\index{smooth manifold} $\M$ and $m\ \in\ \M\,$. Then, the following diagram commutes\index{commutative diagram}:
\begin{equation*}
  \begin{tikzcd}[row sep=large, column sep=large]
    \g \arrow[r, "\mu_{\ast}"] \arrow[d, "\exp"'] & \X\,(\,\M\,) \arrow[d, "\exp"] \\
    G \arrow[r, "\mu_{m}"'] & \M
  \end{tikzcd}
\end{equation*}
or, in other words:
\begin{equation*}
  \forall\,\Xf\ \in\ \g\,: \qquad \mu_{m}(\,\ue^{\,t\,\Xf}\,)\ =\ \ue^{\,t\,\mu_{\ast}(\,\Xf\,)_{m}}\,.
\end{equation*}
\end{proposition}

\subsection{Lie group and Lie algebra representations}
\index{Lie group}\index{Lie algebra}\index{representation}

Let $G$ be a \nm{Lie} group\index{Lie group}, $\g\ =\ \Lie\,(G)$ be its \nm{Lie} algebra\index{Lie algebra} and $V$ be a real vector space\index{vector space}.

\begin{definition}
A \emph{representation}\index{representation} of the \nm{Lie} group\index{Lie group} $G$ in the vector space\index{vector space} $V$ is a homomorphism\index{homomorphism} of \nm{Lie} groups\index{Lie group} (\ie a smooth group\index{smooth group} morphism\index{group morphism}) $\varphi\,:\ G\ \longrightarrow\ \GL\,(V)\,$.
\end{definition}

\begin{definition}
A representation\index{representation} of the \nm{Lie} algebra\index{Lie algebra} $\g$ in the vector space\index{vector space} $V$ is a \nm{Lie} algebra\index{Lie algebra} homomorphism\index{homomorphism} $\phi\,:\ \g\ \longrightarrow\ \End\,(V)\,$.
\end{definition}
Here $\End\,(V)\,$ is enabled with the \nm{Lie} algebra\index{Lie algebra} structure given by the endomorphisms\index{endomorphism} commutator\index{commutator}:
\begin{equation*}
  \forall\, A,\,B\ \in\ \End\,(V)\:\ \qquad \lb{A}{B}\ \eqdef\ A\cdot B\ -\ B\cdot A\,.
\end{equation*}

If $\varphi\,:\ G\ \longrightarrow\ \GL\,(V)$ is a \nm{Lie} group\index{Lie group} representation\index{representation}, then
\begin{equation*}
  \phi\ \eqdef\ \T_{e}\varphi\,:\ \T_{e}G\ =\ \g\ \longrightarrow\ \T_{\mathrm{id}}\bigl(\,\GL\,(V)\,\bigr)\ =\ \End\,(V)
\end{equation*}
is a representation\index{representation} of the \nm{Lie} algebra\index{Lie algebra} $\g\,$.

\subsubsection{Adjoint representations}
\index{adjoint representation}

Let $g\ \in\ G\,$, $V\ =\ \g\,$, then the composition\index{composition} $L_{g}\comp R_{g^{-1}}\,:\ G\ \longrightarrow\ G$ induces a linear mapping\index{linear mapping} $\T_{e}\,(\,L_{g}\comp R_{g^{-1}}\,)\ \eqdef\ \Ad\,(\,g\,)\,:\ \g\ \longrightarrow\ \g$ and, hence, a group morphism\index{group morphism}\index{morphism} $\Ad\,:\ G\ \longrightarrow\ \GL\,(\,\g\,)\,$. Then, the following Lemma holds:

\begin{lemma}
The group morphism\index{group morphism} $\Ad$ is a smooth map\index{smooth map} which gives a representation of $G$ in $\g\,$, which is called the \emph{adjoint} \nm{Lie} group\index{Lie group} representation\index{adjoint Lie group representation}.
\end{lemma}
\begin{proof}
Left to the reader as an exercise.
\end{proof}

Let $\ad\ \eqdef\ \T_{e}\,(\,\Ad\,)\,:\ \g\ \longrightarrow\ \End(\,\g\,)\,$. Then, $\ad$ is also called the \emph{adjoint} \nm{Lie} group\index{Lie group} representation\index{adjoint Lie group representation}.

\begin{lemma}
\begin{equation*}
  \ad\,(\,\Xf\,)\,(\,\Yf\,)\ =\ \lb{\Xf}{\Yf}\,, \qquad \forall \Xf,\,\Yf\ \in\ \g\,.
\end{equation*}
\end{lemma}
\begin{proof}
Left to the reader as an exercise.
\end{proof}

\subsubsection{Co-adjoint representations}
\index{co-adjoint representation}

Let $g\ \in\ G\,$, $V\ \in\ \g^{\ast}$ and $f^{\ast}\ \in\ \End\,(\,\g^{\ast}\,)$ is defined by $f^{\ast}\,(\xi)\ \eqdef\ \xi\comp f$ for any element $f\ \in\ \End\,(\,\g\,)\,$. Then, we can write down the following
\begin{definition}
The following smooth map\index{smooth map}
\begin{align*}
  \Ad^{\,\ast}\,:\ G\ &\longrightarrow\ \GL\,(\,\g^{\,\ast}\,)\,, \\
  g\ &\mapsto\ \Ad\,(\,g^{-1}\,)^{\,\ast}\,,
\end{align*}
which gives a representation of $G$ in $\g^{\ast}$ is called the co-adjoint \nm{Lie} group representation\index{co-adjoint Lie group representation}.
\end{definition}

The last definition makes sense because $\Ad^{\,\ast}\ =\ F\comp\Ad$ and $F\,(f)\ =\ f^{\,\ast}\,$, where the map\index{map} $F\,:\ \End\,(\,V\,)\ \longrightarrow\ \End\,(\,V^{\ast}\,)\,$. Similarly, we can also define
\begin{align*}
  \ad^{\,\ast}\,:\ \g\ &\longrightarrow\ \End\,(\,\g^{\ast}\,)\,, \\
  \Xf\ &\mapsto\ -\ad^{\,\ast}\,(\,\Xf\,)\,,
\end{align*}
where
\begin{equation*}
  \ad^{\,\ast}\,(\,\Xf\,)\,(\,\xi\,)\,(\,\Yf\,)\ =\ -\,\Inner*{\xi}{\lb{\Xf}{\Yf}}\,.
\end{equation*}
Then, $\ad^{\,\ast}$ is also called the co-adjoint \nm{Lie} algebra\index{Lie algebra} representation\index{co-adjoint Lie group representation} in $\g^{\ast}\,$.

\begin{lemma}
\begin{equation*}
  \lb{\ad^{\:\ast}\,(\,\Xf\,)}{\ad^{\:\ast}\,(\,\Yf\,)}\ =\ \ad^{\:\ast}\,(\,\lb{\Xf}{\Yf}\,)\,, \qquad \forall \Xf\,, \Yf\ \in\ \g\,, \quad \forall\xi\ \in\ \g^{\ast}\,.
\end{equation*}
\end{lemma}
\begin{proof}
Left to the reader as an exercise.
\end{proof}

\begin{proposition}
The co-adjoint representation\index{co-adjoint representation} $\Ad^{\:\ast}$ of a \nm{Lie} group $G$ gives rise to a co-adjoint left action\index{co-adjoint left action} of $G$ on $\g^{\ast}\,$:
\begin{align*}
  G\times\g^{\ast}\ &\longrightarrow\ \g^{\ast}\,, \\
  (g,\,\xi)\ &\mapsto\ \Ad^{\:\ast}\,(\,\xi\,)\,.
\end{align*}
\end{proposition}

Let $\Xf\ \in\ \g$ and $F_{\Xf}\ \in\ C^{\infty}\,(\,\g^{\ast}\,)$ be the evaluation function\index{evaluation function} defined by $F_{\Xf}\,(\,\xi\,)\ \eqdef\ \xi\,(\,\Xf\,)\,$. Then, the following Propositions hold:
\begin{proposition}
\begin{equation*}
  \od{}{t} F\,\bigl(\,\Ad^{\:\ast}_{\:\exp\,(\,t\,\Xf\,)}\,(\,\xi\,)\,\bigr)\,\Bigr\vert_{\,t\,=\,0}\ =\ \pb*{F}{F_{\Xf}}\,(\,\xi\,)\,.
\end{equation*}
\end{proposition}

\begin{proposition}
A function $F\ \in\ C^{\infty}\,(\,\g^{\ast}\,)$ is a \nm{Casimir} function\index{Casimir function} for the \nm{Lie}--\nm{Poisson} structure\index{Lie--Poisson structure} on $\g^{\ast}$ if and only if
\begin{equation*}
  \pb*{F}{F_{\Xf}}\ =\ 0\,, \qquad \forall \Xf\ \in\ \g\,.
\end{equation*}
\end{proposition}

Finally, we can state without the proof the following important
\begin{theorem}[\nm{Lie--Berezin--Kirillov--Kostant--Souriau}]
\begin{equation*}
  \Cas\,\bigl(\,C^{\infty}\,(\,\g^{\ast}\,)\,\bigr)\ =\ C^{\infty}\,(\,\g^{\ast}\,)^{\Ad_{\,G}^{\:\ast}}\,,
\end{equation*}
where
\begin{equation*}
  C^{\infty}\,(\,\g^{\ast}\,)^{\Ad_{\,G}^{\:\ast}}\ \eqdef\ \Set*{F\ \in\ C^{\infty}\,(\,\g^{\ast}\,)}{F\comp\Ad_{\,g}^{\:\ast}\ =\ F\,, \quad \forall g\ \in\ G}\,.
\end{equation*}
\end{theorem}

Denote by $\O_{\,\xi}\ \eqdef\ G\cdot\xi\ \subseteq\ \g^{\ast}$ a co-adjoint orbit\index{co-adjoint orbit} of a co-vector $\xi\ \in\ \g^{\ast}\,$. Recall that these orbits are manifolds\index{manifold} such that $\O_{\,\xi}$ admits a submersion\index{submersion} $\phi\,:\ G\ \longrightarrow\ \O_{\,\xi}$ and an immersion\index{immersion} $\imath\,:\ \O_{\,\xi}\ \hookrightarrow\ \g^{\ast}\,$. Define a $\g-$valued $1-$form\index{$1-$form} $\omega$ on $G$ as
\begin{equation*}
  \omega_{\,g}\,(\,\Xf\,)\ =\ \T_{g}\,L_{g^{-1}}\,(\,\Xf\,)\ \in\ \g\,.
\end{equation*}
Then, $L_{\,h}^{\,\ast}\,(\,\omega\,)\ =\ \omega\,$. In other words, the $1-$form\index{$1-$form} $\omega$ is $G-$invariant. Here we take $g,\,h\ \in\ G$ and $\Xf\ \in\ \g\,$. By \nm{Maurer}--\nm{Cartan} formula\index{Maurer--Cartan formula} we have that
\begin{equation*}
  \ud\omega\ =\ -\,\frac{1}{2}\;\lb{\omega}{\omega}\,.
\end{equation*}
Let us define also the $1-$form\index{$1-$form} $\omega_{\,\xi}\ \in\ \Lambda^{1}\,(\,G\,)$ by
\begin{equation*}
  \omega_{\,\xi}\,(\,g\,)\,(\,\Xf\,)\ \eqdef\ \Inner*{\xi}{\omega_{\,g}\,(\,\Xf\,)}\,.
\end{equation*}
Then, the following result can be shown:
\begin{theorem}[\nm{Kirillov--Kostant--Souriau}]
There exists a unique $2-$form\index{$2-$form} $\Omega_{\,\xi}\ \in\ \Lambda^{2}\,(\,\O_{\,\xi}\,)$ such that $\phi^{\ast}\,(\,\Omega_{\,\xi}\,)\ =\ \ud\omega_{\,\xi}\,$. This form is symplectic\index{symplectic form} on the co-adjoint orbit\index{co-adjoint orbit} $\O_{\,\xi}\,$.
\end{theorem}


\section{Moment map, Poisson and Hamiltonian actions}
\label{sec:mm}

\subsection{Introductory motivation}

Let $\R^{3}$ be a basic configuration space\index{configuration space} with coordinate or position vectors\index{position vector} $\r\ =\ (\,q_{1},\,q_{2},\,q_{3}\,)$ and velocity vectors\index{velocity vector}:
\begin{equation*}
  \dr\ =\ (\,\dot{q}_{1},\,\dot{q}_{2},\,\dot{q}_{3}\,)\ \defeq\ \p\ \eqdef\ (\,p_{1},\,p_{2},\,p_{3}\,)\,.
\end{equation*}
We remind that here by the $\dot{(-)}$ we denote the classical time derivative\index{derivative} operator. The angular momentum\index{angular momentum} $\L$ is defined by their vector product\index{vector product} as $\L\ \eqdef\ \r\timesb\dr\,$. The total energy\index{total energy} of a mechanical system is given by
\begin{equation*}
  E_{\mathrm{T}}\ \eqdef\ \frac{\Inner*{\dr}{\dr}}{2}\ +\ U\,(\,\r\,)
\end{equation*}
and the equation of motion is $\ddr\ =\ -\,\grad_{\r}\,U\,(\,\r\,)\,$. The total mechanical energy\index{mechanical energy} is conserved, \ie $\odd{E_{\mathrm{T}}}{t}\ \equiv\ 0\,$. The angular momentum\index{angular momentum} is also constant along a trajectory. It implies that
\begin{equation*}
  \od{\L}{t}\ =\ \dr\timesb\dr\ +\ \r\timesb\ddr\ =\ -\,\r\timesb\grad_{\r}\,U\ =\ \vO\,,
\end{equation*}
which is equivalent to say that $\r\ =\ \lambda\,\grad_{\r}\,U$ for some $\lambda\ \in\ \R\,$.

Let $\so\,(\,3\,)$ be the \nm{Lie} algebra\index{Lie algebra} of skew-symmetric\index{skew-symmetric matrix} $3 \timesb 3$ matrices with real entries. This is a three-dimensional vector space\index{vector space} with the basis\index{basis} $\set*{X_{1},\,X_{2},\,X_{3}}$ given by three following matrices\index{matrix}:
\begin{equation*}
  X_{1}\ =\ \begin{pmatrix}
    0 & 1 & 0 \\
    -1 & 0 & 0 \\
    0 & 0 & 0
  \end{pmatrix}\,, \quad
  X_{2}\ =\ \begin{pmatrix}
    0 & 0 & -1 \\
    0 & 0 & 0 \\
    1 & 0 & 0
  \end{pmatrix}\,, \quad
  X_{3}\ =\ \begin{pmatrix}
    0 & 0 & 0 \\
    0 & 0 & 1 \\
    0 & -1 & 0
  \end{pmatrix}\,.
\end{equation*}
The \nm{Lie} brackets\index{Lie bracket} in the \nm{Lie} algebra\index{Lie algebra} $\so\,(\,3\,)$ are given by
\begin{equation*}
  \lb{X_{i}}{X_{j}}\ =\ X_{k}\,, \qquad (\,i,\,j,\,k\,)\ =\ (\,1,\,2,\,3\,)\,,
\end{equation*}
with all circular permutations\index{permutation}\index{circular permutation}. The \nm{Killing} form\index{Killing form} $\K\,(-,\,-)$ defined as
\begin{align*}
  \K\,:\ \so\,(\,3\,)\times\so\,(\,3\,)\ &\longrightarrow\ \R \\
  (\,X,\,Y\,)\ &\mapsto\ \tr\,(\,X\,Y\,)
\end{align*}
is symmetric\index{symmetric form}, bi-linear\index{bi-linear form} and non-degenerate\index{non-degenerate form}. Here $\tr\,(\,-\,)$ is the trace form of a square matrix. This form identifies $\so\,(\,3\,)$ and $\so\,(\,3\,)^{\ast}$ by the interior product\index{interior product} rule $X\intprod \K\,$.

The \nm{Lie}--\nm{Poisson} structure\index{Lie--Poisson structure} on $\so\,(\,3\,)^{\ast}$ in the coordinates\index{coordinate} $(\,x_{1},\,x_{2},\,x_{3}\,)$ on $\so\,(\,3\,)$ can be expressed as
\begin{equation*}
  \pb*{F}{G}\,(\,x_{1},\,x_{2},\,x_{3}\,)\ =\ \sum_{i,\,j,\,k\, =\, 1}^{3}\,c^{k}_{ij}\,\Bigl(\,\pd{F}{x_{i}}\;\pd{G}{x_{j}}\ -\ \pd{F}{x_{j}}\;\pd{G}{x_{i}}\,\Bigr)\, x_{k}\,.
\end{equation*}
Here $c^{k}_{ij}$ is the structure constant tensor\index{tensor} of the \nm{Lie} algebra\index{Lie algebra} $\so\,(\,3\,).$

The angular momentum\index{angular momentum} $\L$ s defined as a map
\begin{align*}
  \T^{\ast}\R^{3}\ \simeq\ \R^{6}\ &\longrightarrow\ \so\,(\,3\,)^{\ast} \\
  (\,\q,\,\p\,)\ &\mapsto\ \q\timesb\p\ =\ \sum_{i,\,j,\,k}\, (\,q_{\,i}\,p_{j}\ -\ p_{\,i}\,q_{j}\,)\,X_{k}\,.
\end{align*}
The angular momentum\index{angular momentum} map $\L\,:\ \T^{\ast}\R^{3}\ \longrightarrow\ \so\,(\,3\,)^{\ast}$ is a \nm{Poisson} morphism\index{Poisson morphism}.

\begin{definition}
Let $\mu\,:\ G\times\M\ \longrightarrow\ \M$ be a \nm{Lie} group action\index{Lie group action} on a \nm{Poisson} manifold\index{Poisson manifold} $\bigl(\,\M;\,\pb*{}{}\,\bigr)\,$. This action\index{group action} is called a \nm{Poisson} action\index{Poisson action} if the map\index{map}
\begin{equation*}
  \mu_{\,g}^{\ast}\,:\ C^{\infty}\,(\,\M,\,\R\,)\ \longrightarrow\ C^{\infty}\,(\,\M,\,\R\,)
\end{equation*}
defined by $\mu_{\,g}^{\ast}\,(\,F\,)\,(\,m\,)\ \eqdef\ F\,\bigl(\,\mu_{\,g}\,(\,m\,)\,\bigr)$ satisfies the following condition:
\begin{equation*}
  \mu_{\,g}^{\ast}\,\bigl(\,\pb*{F}{G}\,\bigr)\,(\,m\,)\ =\ \pb*{\mu_{\,g}^{\ast}\,(\,F\,)}{\mu_{\,g}^{\ast}\,(\,G\,)}\,(\,m\,)\,, \qquad \forall F,\,G\ \in\ C^{\infty}\,(\,\M,\,\R\,)\,.
\end{equation*}
\end{definition}
Let a \nm{Poisson} structure\index{Poisson structure} $(\,\M,\,\pi\,)$ be symplectic\footnote{Here we mean that the bi-vector $\pi$ is non-degenerate, \ie $\pi$ is invertible when it is seen as a banal matrix.}\index{symplectic structure}. In this case this \nm{Poisson} action\index{Poisson action} can be called the \nm{Hamiltonian} action\index{Hamiltonian action}.

\subsection{Momentum map}

\begin{definition}
Let $\g$ be a \nm{Lie} algebra\index{Lie algebra} and $\bigl(\,\M,\,\pb*{}{}\,\bigr)$ be a \nm{Poisson} manifold\index{Poisson manifold}. A \emph{momentum map}\index{momentum map} is a \nm{Poisson} morphism\index{Poisson morphism} $\mu\,:\ \M\ \longrightarrow\ \g^{\ast}\,$. In other words, it is a smooth map\index{smooth map} $\mu$ such that for $\forall\, F,\,G\ \in\ C^{\infty}\,(\,\g^{\ast}\,)\,$:
\begin{equation*}
  \mu^{\ast}\,\bigl(\,\pb*{F}{G}_{\:\g^{\ast}}\,\bigr)\ =\ \pb*{\mu^{\ast}\,(\,F\,)}{\mu^{\ast}\,(\,G\,)}_{\:\M}\,.
\end{equation*}
\end{definition}
Let $\bar{\lambda}\,:\ \g\ \longrightarrow\ C^{\infty}\,(\,\M,\,\R\,)$ be a smooth linear map\index{linear map}. Then, there is a unique map $\lambda\,:\ \M\ \longrightarrow\ \g^{\ast}$ defined by $\bar{\lambda}\,$:
\begin{equation*}
  \pb*{\lambda\,(\,m\,)}{\Xf}\ =\ \bar{\lambda}\,(\,\Xf\,)\,(\,m\,)\,, \qquad \forall\,m\ \in\ \M\,, \quad \forall\,\Xf\ \in\ \g\,.
\end{equation*}

\begin{proposition}
Let $\bigl(\,\M,\,\pb*{}{}\,\bigr)$ be a \nm{Poisson} manifold\index{Poisson manifold} and $\mu\,:\ \M\ \longrightarrow\ \g^{\ast}$ is a smooth map\index{smooth map}. Then, $\mu$ is a momentum map\index{momentum map} if and only if the associated map $\bar{\mu}\,:\ \g\ \longrightarrow\ C^{\infty}\,(\,\M,\,\R\,)$ is a \nm{Lie} algebra homomorphism\index{Lie algebra homomorphism}:
\begin{equation*}
  \bar{\mu}\,\bigl(\,\lb{\Xf}{\Yf}\,\bigr)\ =\ \pb*{\bar{\mu}\,(\,\Xf\,)}{\bar{\mu}\,(\,\Yf\,)}_{\;\M}\,, \qquad \forall\,\Xf,\,\Yf\ \in\ \g\,.
\end{equation*}
\end{proposition}
Recall that the map $\chi\,:\ C^{\infty}\,(\,\M,\,\R\,)\ \longrightarrow\ \X\,(\,\M\,)$ such that
\begin{equation*}
  \chi\,(\,F\,)\ =\ \Xf_{F}\ =\ \pb*{F}{-}
\end{equation*}
is a \nm{Lie} algebra morphism\index{Lie algebra morphism}. We take the composition $\Theta\ \eqdef\ \chi\comp\bar{\mu}\,:\ \g\ \longrightarrow\ \X\,(\,\M\,)\,$, where $\bar{\mu}\,:\ \g\ \longrightarrow\ C^{\infty}\,(\,\M,\,\R\,)$ and $\g\ =\ \Lie\,(\,G\,)$ with a simply connected \nm{Lie} group\index{Lie group} $G\,$. For compact manifolds\index{compact manifold} $\M\,$, \cref{thm:palais} ensures the existence of an action\index{group action} $\lambda\,:\ G\times\M\ \longrightarrow\ \M$ with $\lambda_{\ast}\ =\ -\,\Theta\,$.

\begin{proposition}
If $G$ is connected, then $\lambda_{\ast}\ =\ -\,\Theta$ gives a \nm{Poisson} morphism\index{Poisson morphism} $\lambda_{\,g}^{\ast}\,:\ C^{\infty}\,(\,\M,\,\R\,)\ \longrightarrow\ C^{\infty}\,(\,\M,\,\R\,)$ for $\forall\,g\ \in\ G$ and $\forall\,u,\,v\ \in\ C^{\infty}\,(\,\M,\,\R\,)\,$:
\begin{equation*}
  \pb*{\lambda_{g}^{\ast}\,(\,u\,)}{\lambda_{g}^{\ast}\,(\,v\,)}\ =\ \lambda_{g}^{\ast}\bigl(\,\pb*{u}{v}\,\bigr)\,.
\end{equation*}
\end{proposition}

\begin{proposition}
Let $\M$ be a compact manifold\index{compact manifold} and $G$ is connected\index{connected group} and simply connected\index{simply connected group}. Then, the action\index{group action} $\lambda$ is $G-$equivariant:
\begin{equation*}
  \begin{tikzcd}[row sep=large, column sep=large]
    \M \arrow[r, "\mu"] \arrow[d, "\lambda_{g}"'] & \g^{\ast} \arrow[d, "\Ad_{\,\g}^{\,\ast}"] \\
    \M \arrow[r, "\mu"'] & \g^{\ast}
  \end{tikzcd}
\end{equation*}
\end{proposition}

\subsection{Moment and Hamiltonian actions}
\index{moment}\index{Hamiltonian action}

Let $\bigl(\,\M,\,\omega\,\bigr)$ be a symplectic manifold\index{symplectic manifold} and the corresponding \nm{Poisson} brackets\index{Poisson bracket} are defined by a pair of \nm{Hamiltonian} vector fields\index{Hamiltonian vector field}:
\begin{equation*}
  \pb*{u}{v}\ =\ \Xf_{\,u}\,(\,v\,)\ =\ \omega\,\bigl(\,\Xf_{\,v},\,\Xf_{\,u}\,\bigr)\,, \qquad \Xf_{\,u}\intprod\omega\ =\ \ud u\,.
\end{equation*}

\begin{lemma}
If $H^{1}\,(\,\M,\,\R\,)\ =\ 0$ and $\Xf\ \in\ \X\,(\,\M\,)$ ``infinitesimally'' preserves the symplectic form\index{symplectic form} $\omega\,$, \ie $\mathcal{L}_{\Xf}\,(\,\omega\,)\ =\ 0\,$, then there exists a unique $u\ \in\ C^{\infty}\,(\,\M,\,\R\,)$ such that $\Xf\ =\ \Xf_{\,u}\,$.
\end{lemma}
Here we should remark that the function $u$ is uniquely defined only modulo a locally constant function on $M$ (which is usually identified with an element of $H^{0}\,(\,\M,\,\R\,).$
\begin{lemma}
Let $\lambda\,:\ G\times\M\ \longrightarrow\ \M$ be an action\index{group action} of a \nm{Lie} group\index{Lie group} $G$ on a symplectic manifold\index{symplectic manifold} $\bigl(\,\M,\,\omega\,\bigr)\,$. The action $\lambda$ is a \nm{Poisson} (more precisely, in this case we may call it a \nm{Hamiltonian}\index{Hamiltonian action}) action\index{Poisson action} if and only if $\lambda_{g}^{\ast}\ =\ \omega\,$.
\end{lemma}

\begin{proposition}
Let $\lambda\,:\ G\times\M\ \longrightarrow\ \M$ be a \nm{Hamiltonian} action\index{Hamiltonian action} on a symplectic manifold\index{symplectic manifold} $\bigl(\,\M,\,\omega\,\bigr)$ and $\lambda_{\ast}\,:\ \g\ \longrightarrow\ \X\,(\,M\,)$ is the corresponding \nm{Lie} algebra homomorphism\index{Lie algebra homomorphism}. Then, $\forall\,\Xf\ \in\ \g\,$:
\begin{equation*}
  \mathcal{L}_{\,\lambda_{\ast}\,(\,\Xf\,)}\,(\,\omega\,)\ =\ 0\,.
\end{equation*}
\end{proposition}

\begin{definition}
Let $\lambda\,:\ G\times\M\ \longrightarrow\ \M$ be an action\index{group action} of $G$ on $\M$ and $\bigl(\,\T^{\ast}\M,\,\Omega\,\bigr)$ is the co-tangent bundle\index{co-tangent bundle} with the canonical symplectic form\index{symplectic form} $\Omega\ =\ \ud\rho\,$, where $\rho$ is the \nm{Liouville} $1-$form\index{Liouville form}\index{$1-$form}. This action\index{group action} can be lifted\index{lift} to an action\index{action} $\Lambda\,:\ G\times\T^{\ast}\M\ \longrightarrow\ \T^{\ast}\M$ defined by
\begin{equation*}
  \Lambda\,(\,g,\,\xi_{\,m}\,)\ \eqdef\ (\,\T_{g\cdot m}\,\lambda_{g^{-1}}^{\ast}\,)\,(\,\xi_{\,m}\,)\,.
\end{equation*}
\end{definition}

\begin{theorem}
The action\index{group action} $\Lambda$ is \nm{Hamiltonian}\index{Hamiltonian action} and the induced momentum map\index{momentum map} $\mu_{\Lambda}\,:\ \T^{\ast}\M\ \longrightarrow\ \g^{\ast}$ is defined by
\begin{equation*}
  \pb*{\mu_{\Lambda}\,(\,\xi_{\,m}\,)}{\Xf}\ =\ \pb*{\xi_{\,m}}{\T_{\,e}\,\lambda_{\,m}\,(\,\Xf\,)}\,.
\end{equation*}
\end{theorem}

\subsubsection{Examples}

\paragraph*{Example 1.} Lifting\index{lifting} of the left\index{left action} $G-$action \emph{on} $G$ to $\T^{\ast}G\,$:
\begin{align*}
  \lambda\,:\ G\times G\ &\longrightarrow\ G\,, \\
  (\,g,\,h\,)\ &\mapsto\ g\cdot h\,.
\end{align*}
Then, we obtain the required lifting\index{lifting}:
\begin{align*}
  \Lambda\,:\ G\times\T^{\ast}G\ \simeq\ G\times\g^{\ast}\ &\longrightarrow\ \T^{\ast}G\ \simeq\ G\times\g^{\ast}\,, \\
  \bigl(\,g,\,(\,h,\,\xi\,)\,\bigr)\ &\mapsto\ (\,g\cdot h,\,\xi\,)\,.
\end{align*}
The associated momentum\index{momentum} can be also easily computed:
\begin{equation*}
  \mu\,(\,\xi_{\,h}\,)\ =\ -\,(\,\T_{e} R_{\,h}\,)^{\ast}(\,\xi_{\,h}\,)\,, \qquad \mu\,(\,h,\,\xi\,)\ =\ \Ad_{\,h}^{\,\ast}\,(\,\xi\,)\,.
\end{equation*}

Similarly, we can consider lifting\index{lifting} of the right\index{right action} $G-$action \emph{on} $G$ to $\T^{\ast}G\,$:
\begin{align*}
  \lambda\,:\ G\times G\ &\longrightarrow\ G\,, \\
  (\,g,\,h\,)\ &\mapsto\ h\cdot g^{-1}\,.
\end{align*}
Then, we obtain the required lifting\index{lifting}:
\begin{align*}
  \Lambda\,:\ G\times\T^{\ast}G\ \simeq\ G\times\g^{\ast}\ &\longrightarrow\ \T^{\ast}G\ \simeq\ G\times\g^{\ast}\,, \\
  \bigl(\,g,\,(\,h,\,\xi\,)\,\bigr)\ &\mapsto\ (\,h\cdot g^{-1},\,\Ad_{\,g}^{\,\ast}\,\xi\,)\,.
\end{align*}
The associated momentum\index{momentum} can be also easily computed:
\begin{equation*}
  \mu\,(\,\xi_{\,h}\,)\ =\ -\,(\,\T_{\,e} L_{\,h}\,)^{\ast}(\,\xi_{\,h}\,)\,, \qquad \mu\,(\,h,\,\xi\,)\ =\ -\,\xi\,.
\end{equation*}

\paragraph*{Example 2.} Let us consider also the action\index{action} of $\S^{1}$ on $\C\,$:
\begin{equation*}
  \C\ \simeq\ \R^{2}\ \simeq\ \T^{\ast}\R^{1}\,, \qquad \Omega\ =\ \ud q\wedge\ud p\,.
\end{equation*}
The action\index{action} is given by
\begin{align*}
  \lambda\,:\ \S^{1}\times\C\ &\longrightarrow\ \C\,, \\
  \bigl(\,\ue^{\,\ui\,\theta},\,z\,\bigr)\ &\mapsto\ \ue^{\,\ui\,\theta}\,z\,,
\end{align*}
for some $\theta\ \in\ [\,0,\,2\,\pi[\,$. Above $\ui$ is the complex imaginary unit, \ie $\ui^{2}\ =\ -\,1\,$. Then, one can easily obtain the expression for $\lambda_{\ast}\,:\ \T_{e}\S^{1}\ \longrightarrow\ \C\,$:
\begin{equation*}
  \lambda_{\ast}\,\Bigl(\,\od{}{\theta}\,\Bigr)\,(q,\,p)\ =\ \T_{e}\,(\,\lambda_{q,\,p}\,)\,\Bigl(\,\od{}{\theta}\,\Bigr)\ =\ -\,p\;\pd{}{q}\ +\ q\;\pd{}{p}\,.
\end{equation*}
The interior product\index{interior product} with the symplectic form\index{symplectic form} can be also easily obtained:
\begin{equation*}
  \lambda_{\ast}\,\Bigl(\,\od{}{\theta}\,\Bigr)\intprod\Omega\ =\ -\,(\,p\ud p\ +\ q\ud q\,)\ =\ -\,\frac{1}{2}\;\ud(\,p^{2}\ +\ q^{2}\,)\,.
\end{equation*}
The momentum map\index{momentum map} is given by
\begin{equation*}
  \mu\,(\,z\,)\ =\ \mu\,(\,q\ +\ \ui\,p\,)\ =\ \frac{p^{2}\ +\ q^{2}}{2}\,.
\end{equation*}
The last construction can be easily generalized to $\C^{\,n}\,$:
\begin{align*}
  \lambda\,:\ \S^{1}\times\C^{\,n}\ &\longrightarrow\ \C^{\,n}\,, \\
  \bigl(\,\ue^{\,\ui\,\theta},\,(\,z_{1},\,z_{2},\ldots,\,z_{n}\,)\,\bigr)\ &\mapsto\ \bigl(\,\ue^{\,\ui\,\theta} z_{1},\,\ue^{\,\ui\,\theta} z_{2},\ldots,\,\ue^{\,\ui\,\theta} z_{n}\,\bigr)\,.
\end{align*}
Then, the associated momentum map\index{momentum map} is given by:
\begin{equation*}
  \mu\,(\,z_{1},\,z_{2},\ldots,\,z_{n}\,)\ =\ \sum_{i\,=\,1}^{n}\,\abs{z_{i}}^{2}\,.
\end{equation*}

\paragraph*{Example 3.} Let us consider the action\index{group action} of $\S^{1}$ on $\S^{2}\,$. The manifold\index{manifold} $\S^{2}$ is equipped with local coordinates\index{local coordinate} $(\,z,\,\phi\,)$ and $\Omega\ =\ \ud z\wedge\ud\phi\,$. The action\index{group action} of $\S^{1}$ is given by rotation\index{rotation} in $z-$planes:
\begin{align*}
  \lambda\,:\ \S^{1}\times\S^{2}\ &\longrightarrow\ \S^{2}\,, \\
  \bigl(\,\ue^{\,\ui\,\theta},\,(\,z,\,\phi\,)\,\bigr)\ &\mapsto\ (\,z,\,\phi\ +\ \theta\,)\,.
\end{align*}
It is not difficult to see that
\begin{equation*}
  \lambda_{\ast}\,\Bigl(\,\od{}{\theta}\,\Bigr)\,(\,z,\,\phi\,)\ =\ \pd{}{\phi}\,, \qquad \lambda_{\ast}\,\Bigl(\,\od{}{\theta}\,\Bigr)\intprod\Omega\ =\ \ud z\,.
\end{equation*}
Finally, the momentum map\index{momentum map} is
\begin{equation*}
  \mu\,(\,z,\,\phi\,)\ =\ z\,.
\end{equation*}

\paragraph*{Example 4.} We consider now the action\index{group action} of $\S^{1}$ on the torus\index{torus} $\T^{2}\ \eqdef\ \S^{1}\times\S^{1}\,$. The torus\index{torus} $\T^{2}$ is equipped with local coordinates\index{local coordinate} $(\,\phi_{1},\,\phi_{2}\,)$ and the symplectic form\index{symplectic form} is $\Omega\ =\ \ud\phi_{1}\wedge\ud\phi_{2}\,$. The action\index{group action} is defined as
\begin{align*}
  \lambda\,:\ \S^{1}\times\T^{2}\ &\longrightarrow\ \T^{2}\,, \\
  \Bigl(\,\ue^{\,\ui\,\theta},\,\bigl(\,\ue^{\,\ui\,\phi_{1}},\,\ue^{\,\ui\,\phi_{2}}\,\bigr)\,\Bigr)\ &\mapsto\ \Bigl(\,\ue^{\,\ui\,\phi_{1}},\,\ue^{\,\ui\,(\,\phi_{1}\ +\ \theta\,)}\,\Bigr)\,.
\end{align*}
Then, we have
\begin{equation*}
  \lambda_{\ast}\,\Bigl(\,\od{}{\theta}\,\Bigr)\,(\,\Omega\,)\ =\ \ud\phi_{1}\wedge\ud(\,\theta\ +\ \phi_{2}\,)\ =\ \Omega
\end{equation*}
and
\begin{equation*}
  \lambda_{\ast}\,\Bigl(\,\od{}{\theta}\,\Bigr)\intprod\Omega\ =\ -\,\ud\phi_{1}\,.
\end{equation*}
Since the coordinate function\index{coordinate function} $\phi_{1}$ is defined only locally, the momentum map\index{momentum map} $\mu$ and the morphism\index{morphism} $\bar{\mu}$ do not exist.

\paragraph*{Example 5.} In this example we consider the action\index{action} of $\SU\,(\,n\,)$ on $\T^{\ast}\bigl(\,\su\,(\,n\,)\,\bigr)\,$. We remind that $\SU\,(\,n\,)$ is the \nm{Lie} group\index{Lie group} of special unitary matrices\index{unitary matrix} with complex coefficients:
\begin{equation*}
  \SU\,(\,n\,)\ \eqdef\ \Set*{\A\ \in\ \Mat_{\,n}(\,\C\,)}{\A\A^{\ast}\ =\ \Id\,,\ \det\,(\A)\ =\ 1}\,,
\end{equation*}
where $\Id$ is the identity matrix\index{identity matrix} and $\A^{\ast}$ is the conjugate (or \nm{Hermitian}\index{Hermitian transpose}) transpose\index{conjugate transpose} of $\A\,$. The corresponding \nm{Lie} algebra\index{Lie algebra} is defined as
\begin{equation*}
  \Lie\,\bigl(\,\SU\,(\,n\,)\,\bigr)\ =\ \su\,(\,n\,)\ \eqdef\ \Set{\A\ \in\ \Mat_{\,n}(\,\C\,)}{\A^{\ast}\ =\ -\,\A\,,\ \tr\,(\,\A\,)\ =\ 0}\,.
\end{equation*}
The \nm{Lie} algebra\index{Lie algebra} $\su\,(\,n\,)$ is an example of a semi-simple \nm{Lie} algebra\index{semi-simple Lie algebra} with a \nm{Killing} form\index{Killing form} $\K\,(\,X,\,Y\,)\ =\ 2\,n\tr\,(\,X\,Y\,)$ and $\su\,(\,n\,)\ \simeq\ \su\,(\,n\,)^{\ast}\,$. The action\index{group action} is defined as
\begin{align*}
  \lambda\,:\ \SU\,(\,n\,)\times\T^{\ast}\bigl(\,\su\,(\,n\,)\,\bigr)\ &\longrightarrow\ \T^{\ast}\bigl(\,\su\,(\,n\,)\,\bigr)\,, \\
  \bigl(\,g,\,(\,X,\,L\,)\,\bigr)\ &\mapsto\ \bigl(\,g\,X\,g^{-1},\,g\,L\,g^{-1}\,\bigr)\,.
\end{align*}
Then, $\Omega\ =\ \tr\,(\,\ud X\wedge\ud L)$ and the momentum map\index{momentum map} is given by
\begin{equation*}
  \mu\,(\,X,\,L\,)\ =\ \lb{X}{L}\,.
\end{equation*}


\section{Reduction of the phase space}
\label{sec:rp}
\index{phase space reduction}

Let $(\,\M,\,\omega\,)$ be a symplectic manifold\index{symplectic manifold} and $\lambda\,:\ G\times\M\ \longrightarrow\ \M$ is a \nm{Hamiltonian} action\index{Hamiltonian action}, \ie
\begin{equation*}
  \lambda_{\,g}^{\,\ast}\,(\,\omega\,)\ =\ \omega\,, \qquad \forall\,g\ \in\ G\,.
\end{equation*}
We justify the terminology by the following observation:
\begin{lemma}
Assume that there exists a momentum map\index{momentum map} $\mu\,:\ \M\ \longrightarrow\ \g^{\ast}\,$, one necessarily obtains that
\begin{equation*}
  \lambda^{\ast}\,(\,\Yf\,)\ =\ -\,\Xf_{\,\bar{\mu}\,(\,\Yf\,)}\,, \qquad \forall\,\Yf\ \in\ \g\,.
\end{equation*}
\end{lemma}

\begin{definition}
An element $c\ \in\ \g^{\ast}$ is called a \emph{regular}\index{regular element} if $\M_{\,c}\ \eqdef\ \mu^{\,-1}\,(\,c\,)$ is a sub-manifold\index{sub-manifold} in $\M$ and if
\begin{equation*}
  \ker\,(\,\T_{m}\,\mu\,)\ =\ \T_{m}\,\M_{\,c}\,, \qquad \forall\,c\ \in\ \M_{\,c}\,.
\end{equation*}
\end{definition}

\begin{lemma}
Let $G_{\,c}\ \eqdef\ \Set*{g\ \in\ G}{\Ad_{\,g}^{\:\ast}\,(\,c\,)\ =\ c}\,$. If $G$ is connected\index{connected group} and simply connected\index{simply connected group}, then $\forall\,m\ \in\ \M_{\,c}$ and $\forall\,g\ \in\ G_{\,c}\,$:
\begin{equation*}
  g\cdot m\ \in\ \M_{\,c}\,.
\end{equation*}
\end{lemma}
\begin{proof}
Left to the reader as an exercise.
\end{proof}

\begin{remark}
In the case when $c$ is a regular element\index{regular element} and $G$ is connected\index{connected group} and simply connected\index{simply connected group}, the action\index{group action} of $G$ on $\M$ induces an action\index{group action} of the \nm{Lie} sub-group\index{Lie sub-group} $G_{\,c}\ \subseteq\ G$ on the sub-manifold\index{sub-manifold} $\M_{\,c}\ \subseteq\ \M\,$.
\end{remark}

\subsection{The main results}

\begin{theorem}
If $G$ is connected\index{connected group} and simply connected\index{simply connected group} and, in addition:
\begin{itemize}
  \item $c$ is a regular element\index{regular element};
  \item $G_{\,c}$ is compact\index{compact};
  \item $G_{\,c}$ acts on $\M_{\,c}$ by free\index{free action} and transitive\index{transitive action} action\index{group action}.
\end{itemize}
Then, there exists a natural smooth structure\index{smooth structure} on $\M_{\,c}\,/G_{\,c}$ such that the mapping\index{mapping} $\pi_{\,c}\,:\ \M_{\,c}\ \longrightarrow\ \M_{\,c}\,/G_{\,c}$ is a submersion\index{submersion}.
\end{theorem}

\begin{remark}
The quotient space\index{quotient space} $\M_{\,c}\,/G_{\,c}$ is called in this case the \emph{reduced phase space}\index{reduced phase space}.
\end{remark}

\begin{theorem}[\nm{Marsden--Weinstein}]\index{Marsden--Weinstein theorem}
If $G$ is connected\index{connected group} and simply connected\index{simply connected group} and, in addition:
\begin{itemize}
  \item $c$ is a regular element\index{regular element};
  \item $G_{\,c}$ is compact;
  \item $G_{\,c}$ acts on $\M_{\,c}$ by free and transitive action\index{group action}.
\end{itemize}
Then, there exists a unique symplectic $2-$form\index{symplectic form}\index{$2-$form} $\omega_{\,c}$ on $\M_{\,c}\,/G_{\,c}$ such that
\begin{equation*}
  \pi_{\,c}^{\:\ast}\,(\,\omega_{\,c}\,)\ =\ \imath_{\,c}^{\:\ast}\,(\,\omega\,)\,,
\end{equation*}
where $\pi_{\,c}\,:\ \M_{\,c}\ \longrightarrow\ \M_{\,c}\,/G_{\,c}$ is the canonical submersion\index{canonical submersion} and $\imath_{\,c}\,:\ \M_{\,c}\ \longhookrightarrow\ \M$ is the canonical embedding\index{canonical embedding}.
\end{theorem}

The proof of this Theorem is based on the following
\begin{lemma}
Let $m\ \in\ \M_{\,c}\,$. Then, $\T_{m}\M_{\,c}\ =\ (\T_{m}G\cdot m)^{\,\perp}\,$. In other words,
\begin{equation*}
  \T_{m}\M_{\,c}\ =\ \Set*{\Xf\ \in\ \T_{m}\M}{\omega_{\,m}\,(\,\Xf,\,\Yf\,)\ =\ 0\,,\ \forall\,\Yf\ \in\ \T_{m}G\cdot m}\,.
\end{equation*}
\end{lemma}
\begin{proof}
Left to the reader as an exercise.
\end{proof}

\begin{remark}
Observe that $\T_{m}\M_{\,c}\ \bigcap\ \T_{m}G\cdot m\ \neq\ \emptyset\,$. More precisely, $\T_{m}\M_{\,c}\ \bigcap\ \T_{m}G\cdot m\ =\ \T_{m}G_{\,c}\cdot m\,$.
\end{remark}

\begin{corollary}
Let $\Xf_{1},\,\Xf_{2},\,\Yf_{1},\,\Yf_{2}\ \in\ \T_{m}\M_{\,c}\ \subseteq\ \T_{m}\M$ such that
\begin{align*}
  \T_{m}\pi_{\,c}\,(\,\Xf_{1}\,)\ &=\ \T_{m}\pi_{\,c}\,(\,\Xf_{2}\,)\,, \\
  \T_{m}\pi_{\,c}\,(\,\Yf_{1}\,)\ &=\ \T_{m}\pi_{\,c}\,(\,\Yf_{2}\,)\,,
\end{align*}
then, $\omega_{\,m}\,(\,\Xf_{1},\,\Yf_{1}\,)\ =\ \omega_{\,m}\,(\,\Xf_{2},\,\Yf_{2}\,)\,$.
\end{corollary}

\begin{lemma}
Let $m,\,n\ \in\ \M_{\,c}$ such that $\pi_{\,c}\,(\,m\,)\ =\ \pi_{\,c}\,(\,n\,)$ and $\Xf_{1},\,\Xf_{2},\,\Yf_{1},\,\Yf_{2}\ \in\ \T_{m}\M_{\,c}\ \subseteq\ \T_{n}\M_{\,c}\ \subseteq\ \T_{n}\M$ such that
\begin{align*}
  \T_{m}\pi_{\,c}\,(\,\Xf_{1}\,)\ &=\ \T_{n}\pi_{\,c}\,(\,\Xf_{2}\,)\,, \\
  \T_{m}\pi_{\,c}\,(\,\Yf_{1}\,)\ &=\ \T_{n}\pi_{\,c}\,(\,\Yf_{2}\,)\,,
\end{align*}
then, $\omega_{\,m}\,(\,\Xf_{1},\,\Yf_{1}\,)\ =\ \omega_{\,n}\,(\,\Xf_{2},\,\Yf_{2}\,)\,$.
\end{lemma}
\begin{proof}
Left to the reader as an exercise.
\end{proof}

\begin{lemma}
Let $c$ be a regular element\index{regular element} and $\O_{\,c}\ \eqdef\ G\cdot c$ be its co-adjoint orbit\index{co-adjoint orbit}. Then, $\mu^{\,-1}\,(\,\O_{\,c}\,)$ is a sub-manifold\index{sub-manifold} in $\M\,$.
\end{lemma}
\begin{proof}
Left to the reader as an exercise.
\end{proof}

\begin{theorem}
The mapping\index{mapping}
\begin{align*}
  \phi\,:\ \O_{\,c}\ &\longrightarrow\ \M_{\,c}\,/G_{\,c}\,, \\
  m\ &\mapsto\ \pi_{\,c}\,(\,g^{-1}m\,)\,,
\end{align*}
where $\mu\,(\,m\,)\ =\ \Ad_{\,g}^{\:\ast}\,(\,c\,)$ is correctly defined, induces a diffeomorphism\index{diffeomorphism}:
\begin{equation*}
  \Phi\,:\ \pi\,\bigl(\,\mu^{\,-1}\,(\,\O_{\,c}\,)\,\bigr)\ \longrightarrow\ \M_{\,c}\,/G_{\,c}\,.
\end{equation*}
\end{theorem}

\subsection{Example}

In this Section we consider again the action\index{group action} of $\S^{1}$ on $\C^{\,n}\,$, which is defined as
\begin{align*}
  \lambda\,:\ \S^{1}\times\C^{\,n}\ &\longrightarrow\ \C^{\,n}\,, \\
  \bigl(\,\ue^{\,\ui\,\theta},\,\q\ +\ \ui\p\,\bigr)\ &\mapsto\ \ue^{\,\ui\,\theta}\,\q\ +\ \ui\,\ue^{\,\ui\,\theta}\,\p\,,
\end{align*}
where $\p,\,\q\ \in\ \R^{n}\,$. The momentum map\index{momentum map} is
\begin{align*}
  \mu\,:\ \C^{\,n}\ &\longrightarrow\ \Lie\,(\,\S^{1}\,)\,, \\
  \q\ +\ \ui\,\p\ &\mapsto\ -\,\sum_{i\,=\,1}^{n}\,\frac{q_{i}^{\,2}\ +\ p_{i}^{\,2}}{2}\,.
\end{align*}
Then, $\M_{\,c}\ =\ \Set*{z\ \in\ \C^{\,n}}{\sum_{i\,=\,1}^{n}\,\abs{z_{i}}^{\,2}\ =\ 2\,c}\ \simeq\ \S^{\,n}\,$, $c\ >\ 0\,$. It is also clear that $G_{\,c}\ \simeq\ \S^{1}$ and $\S^{1}$ is an \nm{Abelian} group\index{Abelian group}. Thus, we have:
\begin{equation*}
  L_{g}R_{g^{-1}}\ =\ \Id \qquad \Longrightarrow \qquad \Ad_{\,g}\ =\ \Ad_{\,g}^{\:\ast}\ =\ \Id\,.
\end{equation*}
Henceforth,
\begin{equation*}
  \M_{\,c}\,/G_{\,c}\ =\ \S^{\,n}\,/\,\S^{1}\ \simeq\ \P^{\,n\,-\,1}\,.
\end{equation*}


\section{Poisson--Lie groups}
\label{sec:last}

A \nm{Lie} group\index{Lie group} $G$ is called \emph{\nm{Poisson--Lie} group}\index{Poisson--Lie group} if it is a \nm{Poisson} manifold\index{Poisson manifold} such that the multiplication\index{multiplication} $m\,:\ G\times G\ \longrightarrow\ G$ is a morphism\index{morphism} of \nm{Poisson} manifolds\index{Poisson manifold}. Let $\g$ be \nm{Lie} algebra\index{Lie algebra}, $\g^{\ast}$ be dual vector space\index{dual vector space} to $\g\,$. 
 
\begin{definition}
We say that $\g$ is a \emph{\nm{Lie} bi-algebra}\index{Lie bi-algebra} if there is a \nm{Lie} algebra\index{Lie algebra} structure $\lb*{}{}_{\,\ast}$ on $\g^{\ast}$ such that the map\index{mapping} $\delta\,:\ \g\ \longrightarrow\ \Lambda^{2}\,\g$ (called the {\it co-bracket}\index{co-bracket}), dual to the bracket $\lb*{}{}_{\,\ast}\,:\ \Lambda^{2}\,\g^{\ast}\ \longrightarrow\ g^{\ast}$ is a $1-$cocycle with respect to the adjoint action\index{adjoint action} of $\g$ on $\Lambda^{2}\,\g\,$. 
\end{definition}

\subsection{Modified Classical Yang--Baxter equation}
\index{Yang--Baxter equation}

Let $G$ be connected\index{connected group} and a simply connected\index{simply connected group} \nm{Lie} group\index{Lie group}, and let $\g$ be its \nm{Lie} algebra\index{Lie algebra}. Then there is one-to-one correspondence\index{one-to-one correspondence} between \nm{Poisson--Lie} group\index{Poisson--Lie group} structures\index{Poisson--Lie structure} on $G$ and \nm{Lie} bi-algebra\index{Lie bi-algebra} structures on $\g\,$.

As V.~\nm{Drinfel'd} showed \cite{Drinfeld1983}, every \index{Poisson--Lie} structure\index{Poisson--Lie structure} on a semi-simple\index{semi-simple group} connected\index{connected group} $G$ has the following form:
\begin{equation}\label{eq:PL}
  \pi\,(\,g\,)\ =\ \Lambda^{2}\,\Bigl(\,(\,\Ll_{\,g}\,)_{\,\ast}\,\Bigr)\,(\,\r\,)\ -\ \Lambda^{2}\,\Bigl(\,(\,\Rr_{\,g}\,)_{\,\ast}\,\Bigr)\,(\,\r\,)\,,
\end{equation}
where $(\,\Ll_{\,g}\,)_{\,\ast}$ and $(\,\Rr_{\,g}\,)_{\,\ast}$ denote tangent maps\index{tangent map} of left and right translations by $g\ \in\ G\,$. The element\index{element} $\r\ \in\ \Lambda^{2}\,\g$ satisfies the following condition:
\begin{equation}\label{eq:mCYBE}
 \llbracket\,\r,\,\r\,\rrbracket\ \eqdef\ \lb{\r_{\,12}}{\r_{\,13}}\ +\ \lb{\r_{\,12}}{\r_{\,23}}\ +\ \lb{\r_{\,13}}{\r_{\,23}}\ \in\ \Lambda^{3}\,\g\,, 
\end{equation}
where the right hand side is invariant under the adjoint action\index{adjoint action} of $\g\,$. The condition~\eqref{eq:mCYBE} is called a \emph{modified \nm{Yang--Baxter} equation}\index{Yang--Baxter equation} and the bracket\index{bracket}
\begin{equation*}
  \llbracket\,-,\,-\,\rrbracket\,:\ \Lambda^{2}\,\g \otimes \Lambda^{2}\,\g\ \longrightarrow\ \Lambda^{3}\,\g
\end{equation*}
is a so-called \emph{\nm{Schouten--Nijenhuis} bracket}\index{Schouten--Nijenhuis bracket}. This is the natural graded\index{graded structure} (or \emph{super}-) \nm{Lie} algebra\index{Lie algebra} structure on the exterior algebra\index{exterior algebra}
\begin{equation*}
  \Lambda^{\bullet}\,\g\ =\ \bigoplus_{k}\, \Lambda^{k}\,\g\,
\end{equation*}
Here $\r_{\,12}$, to give an example, denotes an element $\r_{\,12}\ =\ \r\: \otimes\: \Id_{\,3}\ \in\ (\,\g\: \otimes\: \k\,)^{\otimes\,3}\,$; $\k\ \in\ \set{\R,\,\C}$ and $\r$ being usually called a \emph{classical $\r-$matrix}.

The condition \eqref{eq:mCYBE} ensures that the bracket\index{Poisson bracket} $\pb{}{}_{\,\ast}$ on $\g^{\ast}$ satisfies the \nm{Jacobi} identity\index{Jacobi identity}. The corresponding \nm{Lie} bi-algebra\index{Lie bi-algebra} structure is calculated in the obvious way. Namely, the co-bracket\index{co-bracket} $\delta$ is given by
\begin{equation*}
  \delta\,(\,x\,)\ =\ \ud_{\,e}\,\pi\,(\,x\,)\ =\ \Ll_{\bar x}\,\pi\,(\,e\,)\ =\ \od{}{t}\;\r_{\,(\,\ue^{\,-t\,x}\,)_{\,\ast}}\,\pi\,(\,\ue^{\,t\,x}\,)\,\Bigr\vert_{\,t\,=\,0}\ =\ \ad_{\,x}\,(\,\r\,)\,,
\end{equation*}
where $\ud_{\,e}\,\pi$ is the intrinsic derivative\index{intrinsic derivative} of a poly-vector field\index{poly-vector field} on $G$ with $\pi\,(\,e\,)\ =\ 0\,$, $\bar x$ is any vector field\index{vector field} on $G$ with $\bar x\,(\,e\,)\ =\ x$ and $\Ll_{\,\bar x}$ denotes the \nm{Lie} derivative\index{Lie derivative} \cite{Lu1990}.

The \nm{Poisson} structures\index{Poisson structure} of the form \eqref{eq:PL} are called \emph{co-boundary}\index{co-boundary structure} or \emph{$\r-$matrix structures}\index{$\r-$matrix structure}. Since for a connected\index{connected group} semi-simple\index{semi-simple group} or a compact\index{compact group} \nm{Lie} group\index{Lie group} $G$ every $1-$cocycle is a co-boundary, one has the following
\begin{proposition}. 
The \nm{Poisson--Lie} structures\index{Poisson--Lie structure} on a connected\index{connected group} semi-simple\index{semi-simple group} or a compact\index{compact group} \nm{Lie} group\index{Lie group} $G$ are in one-to-one correspondence\index{one-to-one correspondence} with the solutions $\r\ \in\ \Lambda^{2}\,\g$ of the modified \nm{Yang--Baxter} equation\index{Yang--Baxter equation}.
\end{proposition}

\subsection{Manin triples}
\index{Manin triple}

Let $\g$ be a \nm{Lie} bi-algebra\index{Lie bi-algebra}. There is a unique \nm{Lie} algebra\index{Lie algebra} structure on the vector space\index{vector space} $\g \oplus \g^{\ast}$ such that
\begin{enumerate}
  \item $\g$ and $\g^{\ast}$ are \nm{Lie} sub-algebras\index{Lie sub-algebra}.
  \item The symmetric\index{symmetric form} bi-linear form\index{bi-linear form} on $\g \oplus \g^{\ast}$ given by the relation
  \begin{equation*}
    \Inner{\Xf\ +\ \xi}{\Yf\ +\ \eta\Yf}\ =\ \Inner{\Xf}{\eta}\ +\ \Inner{\Yf}{\xi}\,, \quad \forall\,\Xf,\,\Yf\ \in\ \g\,, \quad \forall\,\xi,\,\eta\ \in\ \g^{\ast}
  \end{equation*}
  is invariant.
\end{enumerate}
This structure is given by 
\begin{equation*}
  \pb*{\Xf}{\xi}\ =\ -\ad^{\,\ast}_{\,\Xf}\,(\,\xi\,)\ +\ \ad^{\,\ast}_{\,\xi}\,(\,\Xf\,)\,,
\end{equation*}
for $\Xf\ \in\ \g$ and $\xi\ \in\ \g^{\ast}\,$, where $\ad^{\,\ast}$ is the co-adjoint action\index{co-adjoint action}. This \nm{Lie} algebra\index{Lie algebra} is denoted by $\g \bowtie \g^{\ast}$ and $(\,\g \bowtie \g^{\ast},\, \g,\, \g^{\ast}\,)$ is an example of a \nm{Manin} triple\index{Manin triple}. In general, a \nm{Manin} triple\index{Manin triple} is a decomposition\index{Lie algebra decomposition} of a \nm{Lie} algebra\index{Lie algebra} $\g$ with a non-degenerate invariant scalar product\index{scalar product} $\Inner{}{}$ into direct sum of isotropic with respect to $\Inner{}{}$ vector spaces\index{vector space}, $\g\ =\ {\g}_{\,+} \oplus \g_{\,-}$ such that $\g_{\,\pm}$ are \nm{Lie} sub-algebras\index{Lie sub-algebra} of $\g\,$. It is well-known that there is one-to-one correspondence\index{one-to-one correspondence} between \nm{Lie} bi-algebras\index{Lie bi-algebra} and \nm{Manin} triples\index{Manin triple}. These triples were introduced by V.~\nm{Drinfel'd} \cite{Drinfeld1988} and named after Yu.~I.~\nm{Manin}.

\subsection{Poisson--Lie duality}
\index{Poisson--Lie duality}

Let $G$ be a connected\index{connected group} and simply connected\index{simply connected group} \nm{Poisson--Lie} group\index{Poisson--Lie group}, $\g\ =\ \Lie\,(\,G\,)$ its \nm{Lie} algebra\index{Lie algebra} and $\bigl(\,\g\ \bowtie\ \g^{\ast},\, \g,\, \g^{\ast}\,\bigr)$ the \nm{Manin} triple\index{Manin triple}. By duality, $\bigl(\,\g^{\ast}\ \bowtie\ \g,\, \g^{\ast},\, \g\,\bigr)$ is also a \nm{Manin} triple\index{Manin triple}. Then $\g^{\ast}$ is a \nm{Lie} bi-algebra\index{Lie bi-algebra}. This enables us to consider a connected\index{connected group} and simply connected\index{simply connected group} \nm{Lie} group\index{Lie group} $G^{\ast}$ with a \nm{Poisson--Lie} structure\index{Poisson--Lie structure} $\pi^{\ast}$ and with the tangent \nm{Lie} bi-algebra\index{tangent Lie algebra} $\g^{\ast}\,$. The \nm{Poisson--Lie} group\index{Poisson--Lie group} $(\,G^{\ast},\, \pi^{\ast}\,)$ is called the \emph{\nm{Poisson--Lie} dual}\index{Poisson--Lie dual} to $(\,G,\,\pi\,)\,$.


\subsection{Example of non-Hamiltonian action}
\index{non-Hamiltonian action}

Let $G$ be a \nm{Poisson}--\nm{Lie} group\index{Poisson--Lie group} with a multiplicative \nm{Poisson} tensor $\pi_g$ and $\M$ be a smooth \nm{Poisson} manifold\index{Poisson manifold} with a \nm{Poisson} structure\index{Poisson structure} given by $\pi_{\,\M}\,$. Then, the product $G\times\M$ can be considered as a \nm{Poisson} manifold\index{Poisson manifold} with the direct sum\index{direct sum} structure $\tilde{\pi}\,$.

\begin{proposition}
An action $\sigma\,:\ G\times\M\ \longrightarrow\ \M$ of a \nm{Poisson--Lie} group\index{Poisson--Lie group} $G$ on a \nm{Poisson} manifold\index{Poisson manifold} $\M$ is a \nm{Poisson--Lie} action\index{Poisson--Lie action} if and only if
\begin{equation*}
  \pi_{\,\M}\,(\,g\cdot m\,)\ =\ \Lambda^{2}\,\bigl(\,(\,\sigma_{\,g}\,)_{\,\ast}\,\bigr)\,\bigl(\,\pi_{\,\M}\,(\,m\,)\,\bigr)\ +\ \Lambda^{2}\,\bigl(\,(\sigma_{\,m}\,)_{\,\ast}\,\bigr)\,\bigl(\,\pi_{\,G}\,(\,g\,)\,\bigr)\,.
\end{equation*}
\end{proposition}

\begin{remark}
One can consider any \nm{Lie} group\index{Lie group} $G$ as a \nm{Poisson--Lie} group\index{Poisson--Lie group} with $\pi_{\,G}\ \equiv\ 0$ then the action $\sigma$ is a \nm{Poisson} (action)\index{Poisson action} if it gives a \nm{Poisson} morphism\index{Poisson morphism} $\pi_{\,\M}\,(\,g\cdot m\,)\ =\ (\,\sigma_{\,g}\,)_{\,\ast}\,\bigl(\,\pi_{\,\M}\,(\,m\,)\,\bigr)\,$.
\end{remark}

\begin{definition}
The action\index{group action} $\sigma\,:\ G\times\M\ \longrightarrow\ \M$ is called a \nm{Poisson}--\nm{Lie} action\index{Poisson--Lie action} if $\pi^{\ast}\,:\ C^{\infty}\,(\,\M\,)\ \longrightarrow\ C^{\infty}\,(\,G\times\M\,)$ is a \nm{Poisson} morphism\index{Poisson morphism}:
\begin{equation*}
  \pi^{\ast}\,\bigl(\,\pb*{F}{H}_{\;\pi_{\,\M}}\,\bigr)\ =\ \pb*{\pi^{\ast}\,(\,F\,)}{\pi^{\ast}\,(\,H\,)}_{\;\tilde{\pi}}\,.
\end{equation*}
\end{definition}

Infinitesimally, a \nm{Poisson--Lie} action\index{Poisson--Lie action} of a \nm{Lie} bi-algebra\index{Lie bi-algebra} $\g$ on a \nm{Poisson} manifold $(\,\M,\,\pb*{}{}\,)$ is given by an action\index{group action} 
\begin{align*}
  \rho\,:\ \g &\longrightarrow\ \X\,(\,\M\,)\,, \\
  \Xf\ &\mapsto\ V_{\Xf}\,,
\end{align*}
with $\Xf\ \in\ \g$ such that
\begin{multline*}
  V_{\Xf}\,\pb*{f}{g}\,(\,m\,)\ =\ \pb*{V_{\Xf} f}{g}\,(m)\ +\ \pb*{f}{V_{\Xf}\,g}\,(\,m\,)\\ 
  -\ \pb*{\Xf}{\lb*{\rho^{\ast}\,\ud\,f\,(\,m\,)}{\rho^{\ast}\,\ud\,g\,(\,m\,)}}\,,
\end{multline*}
where $\rho^{\ast}\,\ud f\,(\,m\,)\ \in\ \g^{\ast}$ and $\Inner{\Xf}{\rho^{\ast}\,\ud f\,(\,m\,)}\ =\ V_{\Xf} f\,(\,m\,)\,$. In other words,
\begin{equation*}
  \Inner{\Xf}{\lb*{\tilde{\rho}\,(\,\ud F\,)\,(\,m\,)}{\tilde{\rho}\,(\,\ud G\,)\,(\,m\,)}_{\,\ast}}\ =\ \Inner{\ud F}{\ud G}\,\bigl(\,\rho\,(\,\Xf\,)\,\bigr)\,(\,m\,)\,.
\end{equation*}
define a \emph{\nm{Lie} algebroid}\index{Lie algebroid} structure on $\T^{\ast}\M\,$.

There are natural left and right actions\index{group action} of dual \nm{Poisson--Lie} group\index{Poisson--Lie group} $G^{\ast}$ on $G\,$. These actions are called \emph{left (right) dressing} transformations\index{dressing transformation}. The dressing transformations\index{dressing transformation} are not \nm{Hamiltonian}\index{Hamiltonian action} as \nm{Semenov-Tian-Shansky} proved but these actions are genuine \nm{Poisson--Lie} actions\index{Poisson--Lie action} \cite{Semenov-Tian-Shansky1985}.


\subsection*{Acknowledgments}
\addcontentsline{toc}{subsection}{Acknowledgments}

D.D. was supported by the Laboratory of Mathematics (LAMA UMR \#$5127$) and the University \nm{Savoie Mont Blanc} to attend the meeting in \nm{Wis\l{}a}. V.R. acknowledges a partial support of the project IPaDEGAN (H$2020$-MSCA-RISE-$2017$), Grant Number $778010$, and of the Russian Foundation for Basic Research under the Grants RFBR $18-01-00461$ and $16-51-53034-716$ GFEN. Both Authors would like to thank the \nm{Baltic} Mathematical Institute for organizing this scientific event and the anonymous Referee who helped us to improve the presentation indicating some shortcomings and misprints.



\bigskip\bigskip
\invisiblesection{References}
\bibliographystyle{acm}
\bibliography{mybiblio}
\bigskip\bigskip


\newpage
\printindex


\TheEnd


\end{document}